\documentclass[11pt]{article}

\usepackage{amssymb}
\usepackage{amsmath}
\usepackage{amsthm}
\usepackage{epsfig}
\usepackage{mathrsfs}
\usepackage{graphicx}

\def\RR{\mathbb{R}}

\def\DD{{\cal D}}
\def\FF{{\cal F}} 
\def\LL{{\cal L}} 
\def\SS{{\cal S}}

\def\LnF{L^2_{\nu F^{-1}}}

\def\wto{\rightharpoonup}

\def\vphi{\varphi}
\def\pa{\partial}
\def\na{\nabla}
\def\eps{\varepsilon}

\newtheorem{theorem}{Theorem}[section]
\newtheorem{lemma}[theorem]{Lemma}
\newtheorem{proposition}[theorem]{Proposition}

\newtheorem{remark}[theorem]{Remarks}
\newtheorem{rk&ex}[theorem]{Remarks \& Examples}

\def\qed{\hbox{${\vcenter{\vbox{                        
   \hrule height 0.4pt\hbox{\vrule width 0.4pt height 6pt
   \kern5pt\vrule width 0.4pt}\hrule height 0.4pt}}}$}}

\title{Fractional diffusion limit for collisional
  kinetic equations}
\author{Antoine Mellet$^{\mbox{{\footnotesize 1}}}$,
  St\'ephane Mischler$^{\mbox{{\footnotesize 2}}}$ 
  and Cl\'ement Mouhot$^{\mbox{{\footnotesize 2}}}$}

\date{}

\begin{document}
\maketitle

\footnotetext[1]{Email: \texttt{mellet@math.umd.edu} \\ Department of
  Mathematics, University of Maryland, College Park, MD 20742 USA}
\footnotetext[2]{Emails: \texttt{mischler@ceremade.dauphine.fr},
  \texttt{mouhot@ceremade.dauphine.fr} \\
  CEREMADE, Universit\'e Paris Dauphine, Place du Mar\'echal de Lattre
  de Tassigny, 75775 Paris cedex 16, FRANCE}

\vspace{0.3cm}

\begin{abstract}
  This paper is devoted to diffusion limits of linear Boltzmann
  equations.  When the equilibrium distribution function is a
  Maxwellian distribution, it is well known that for an appropriate
  time scale, the small mean free path limit gives rise to a diffusion
  equation.  In this paper, we consider situations in which the
  equilibrium distribution function is a heavy-tailed distribution
  with infinite variance.
  We then show that for an appropriate time scale, the small mean free
  path limit gives rise to a fractional diffusion equation. 
\end{abstract}

\vspace{0.3cm}

Fractional derivatives and integrals.


\textbf{Acknowledgments}: Antoine Mellet gratefully thanks the CEREMADE at the Universit\'e Paris Dauphine, 
where most of this research was performed, for its hospitality.
Antoine Mellet was also partially supported by NSERC Grant 341253-07. 

\tableofcontents


\section{Introduction}


It is well known that under appropriate scaling, the asymptotic
analysis of collisional Êkinetic equations can lead to diffusion-type
equations.  This scaling corresponds to a long time scale and a small
mean-free path.  More precisely, the starting point is the following
collisional kinetic equation:
\begin{eqnarray}
&& \pa_t f + v \cdot \nabla_x f = L(f) \qquad  
    \mbox{ in } (0,\infty) \times \RR^N \times \RR^N, \label{eq:0} \\[4pt]
&&  f (0,.) =f_0 \qquad\qquad\qquad\mbox{ in } \RR^N \times \RR^N ,\label{eq:01}
\end{eqnarray}
which models the evolution of a particle distribution function,
$f(t,x,v) \ge 0$ depending on the time $t>0$, the position $x\in\RR^N$
and a variable $v\in\RR^N$.  This variable $v$ usually represents the
velocity of the particles, or some other internal degree of freedom of
the particles, such as a wave vector.  For simplicity, we take $v\in
\RR^N$, though other spaces could be considered (torus for the wave
vector in semiconductor or bounded set for relativistic particles).

\smallskip 

We then introduce the macroscopic variables
$$ x'=\eps x\qquad  t'= \theta(\eps) \,  t $$ 
and the rescaled distribution function
$$ f^\eps(t',x',v) = f(t,x,v),$$
where $\eps$ is a small parameter. The function $f^\eps$ is now
solution of (we have skipped the primes)
\begin{equation}\label{eq:eps1}
\theta(\eps) \,  \pa_t f^\eps + \eps \, v \cdot \na_x f^\eps =  L(f^\eps).
\end{equation}
The object of this paper is to investigate the behavior, as $\eps$
goes to zero, of the solutions of (\ref{eq:eps1}).  It will naturally
strongly depend on the properties of the collision operator $L$.

\smallskip

Throughout this paper, we will assume that $L$ is a linear Boltzmann
operator (sometimes also called ``scattering'' operator), describing
the interactions of the particles with the surrounding medium, of the
form:
\begin{equation}\label{eq:LL}
L(f) = \int_{\RR^N} \left[ \sigma(v,v') f(v') - \sigma(v',v) f(v)\right]\, dv'
\end{equation}
with a non-negative {\em collision kernel} $\sigma=\sigma(v,v') \ge
0$.  The operator $L$ is conservative, i.e. it preserves the total
mass of the distribution.  Under classical assumptions on the
collision kernel $\sigma$, there exists a unique positive normalized
equilibrium function $F$:
$$F = F(v) > 0 \mbox{ a.e. on $\RR^N$},\quad \int_{\RR^N} F (v)\, dv=1 
\mbox{ and }L(F)=0 .$$ 
In the sequel, we always assume that $F$ exists and is an even function of $v$. 

\smallskip

The derivation of diffusion-type equations from kinetic equations such
as (\ref{eq:eps1}) was first investigated by E. Wigner \cite{W}, A.
Bensoussan, J.L.  Lions and G. Papanicolaou in \cite{BLP} and E.W.
Larsen and J.B. Keller \cite{LK}. In \cite{DGP}, P. Degond, T. Goudon
and F. Poupaud consider very general collision operators of the form
(\ref{eq:LL}).  When $F$ decreases "quickly enough" for large values
of $|v|$ (and under a few additional assumptions on $\sigma$ and $F$),
they prove in particular that for $\theta(\eps) = \eps^2$,
$f^\eps(t,x,v)$ converges, when $\eps$ goes to zero, to a function of
the form $\rho(t,x)F(v)$ where the density $\rho(t,x)$ solves a
diffusion equation
\begin{equation}\label{eq:diffusion}
\pa_t \rho - \na_x \left( D \, \na_x \rho \right) =0,
\end{equation}
with diffusion matrix $D$ given by the following formula
\begin{equation}\label{eq:coefdiffclass}
D= \int_{\RR^N} \left(  v\otimes \chi \right) \, dv\, 
\quad \mbox{  with } \quad L(\chi ) = -v F.
\end{equation}

\smallskip

In order to introduce our problem in simple terms, we now consider, in
this introduction and in the next section, that $L$ has the following
very simple form:
\begin{equation}\label{eq:Lsimple}
 L(f) =\rho \,  F-f,\qquad \rho = \langle f \rangle := \int_{\RR^N} f(v)\, dv,
\end{equation}
which corresponds to the choice of a collision kernel
$\sigma(v,v')=F(v)$. This baby model is usually called ``linear
relaxation'' or sometimes ``linear BGK'' collision operator. We shall
come back to more general (and realistic) collision kernel
$\sigma(v,v')$ in Section~\ref{sectionMR}.  Under (\ref{eq:Lsimple}),
it is readily seen that we can take $\chi=v\, F$ in the definition of
the diffusion matrix $D$ and thus:
$$D= \int_{\RR^N}  \left( v\otimes v \right) \, F(v)\, dv.$$
In particular, in order for $D$ to be finite, we need the second
moment of $F$ to be finite:
\begin{equation}\label{eq:moment2}
 \int_{\RR^N}  |v|^2 F(v)\, dv <\infty.
\end{equation}
This is the case, for example, when the equilibrium $F$ is given by
the so-called ``Maxwellian'' distribution $F(v) =
C\exp(-\frac{v^2}{2})$. The main result of this paper can then be
  summarized as follows: If $F$ is such that (\ref{eq:moment2}) does
  not hold, then, under appropriate time scale, the limit
  $\eps\rightarrow 0$ in (\ref{eq:eps1}) leads to a fractional
  diffusion equation instead of (\ref{eq:diffusion}).

\smallskip

More precisely, we will assume that $F$ is a heavy-tailed distribution
functions, that is typically (we will consider slightly more general
$F$ later on):
\begin{equation}\label{eq:poly} 
F (v) \sim \frac{\kappa_0}{|v|^{N+\alpha}} \quad\hbox{as}\quad |v| \to \infty, 
\end{equation}
with $\kappa_0 > 0$ and $\alpha > 0$ (this last condition guarantees
that $F$ is integrable).

\smallskip

Heavy-tailed distribution functions arise in many contexts. For
instance, most astrophysical plasmas are observed to have velocity
distribution functions exhibiting power law tails (see Summers and
Thorne \cite{ST} or Mendis and Rosenberg \cite{MR}). Dissipative
collision mechanisms in granular gases can also produce power law
tails: see for instance \cite{EB} for the so-called ``inelastic
Maxwell model'' introduced in \cite{BCG}. One can also refers to the
more general review paper \cite{Vill-gran}. We also mention that a
recent work \cite{BG} has shown that even elastic collision mechanisms
can produce power law tail behaviors in the case of mixture of gases
with Maxwellian collision kernel. Power law tails are also common in
economy where they are referred to as Pareto distributions (as well as
in statistics and probability more generally). Here are some samples
of (very different) mathematical works using statistical physics
models to account for these power laws in economy: Newman \cite{N},
Duering-Toscani \cite{DT} and Wright \cite{Wright}. Let us just
mention to the reader that, in the case of economy, the interpretation
of these power law distributions as equilibria of statistical physics
models (which is the core of the field now known as ``econophysics'')
is controversial since it does not take into account any individual
rationality.

It is thus the goal of this paper to investigate what happens to
diffusion limits when the velocity repartition is no longer described
by a Gaussian function. More precisely we shall be particularly
interested in those heavy-tailed distributions whose variance is
infinite, since in this case we shall show that the diffusion limit
yields equations of the form (\ref{eq:D1}), and our underlining motivation was
indeed to provide a microscopic derivation for these fractional diffusion
equations, which is lacking at now.

As a striking link between these goal and motivation, let us mention
that a famous case of such distribution with infinite variance are the
stable (or L\'evy) distributions.  And the latter plays an important role in
probability theory since it can be interpreted as the law of a
``L\'evy walk'' whose law evolution is governed by a fractional diffusion
equation. To say it differently and in a more analytical way, these
stable distributions are the fundamental solutions of the fractional
diffusion equation. They thus play the role played by the gaussian
distribution in the case of the heat equation. 
 
\smallskip

When (\ref{eq:poly}) holds with $\alpha>2$ then (\ref{eq:moment2}) is
still satisfied, and the analysis leading to (\ref{eq:diffusion}) can
be performed. We are thus interested in values of $\alpha$ less than
$2$.  More precisely, when $\alpha\in(0,2)$, we will prove that the
appropriate diffusion scaling is given by
\begin{equation}\label{eq:scal1} 
\theta(\eps) := \eps^\alpha, 
\end{equation}
and that the solution $f^\eps$ of (\ref{eq:eps1}) then converges to
$\rho(t,x) F(v)$ with $\rho$ solution of the following fractional
diffusion equation
\begin{eqnarray}
&& \pa_t \rho + \kappa \, (-\Delta_{x})^{\alpha/2} \rho= 0 \qquad 
\mbox{ in }(0,\infty) \times \RR^N,\label{eq:D1} \\[4pt]
&& \rho(0,.) = \rho_0  \quad \qquad  \qquad\,  
\qquad \mbox{ in }\RR^N,\label{eq:D2}
\end{eqnarray}
with
\begin{equation}\label{eq:D3}
\kappa =  \int_{\RR^N} \frac{ w_1^2}{ 1 + w_1^2 } \, {\kappa_0 \over  |w|^{N+\alpha}} \, dw.
\end{equation}
We recall that the operator $ (-\Delta_{x})^{\alpha/2} $ denotes the
fractional operator, defined for instance by the Fourier formula
$$ 
(-\Delta_{x})^{\alpha/2} \rho := {\cal F}^{-1} \Bigl( |k|^\alpha \, {\cal F}( \rho) (k)\Bigr),
$$  
where ${\cal F}$ stands for the Fourier transform in the space
variable.  This operator has a lot of nice properties, similar to that
of the usual Laplace operator, most notably it retains some
ellipticity (in the sense that its $L^2$ Dirichlet form controls some
fractional Sobolev norm $H^{\alpha/2}$). The fundamental difference is the fact that
this operator is nonlocal for any $0<\alpha<2$.

\smallskip

Before we state our first result, let us introduce some functional
spaces definitions. We denote by $L^p$, $p \in [1,+\infty]$ the usual
Lebesgue space on $\RR^N$ or $\RR^N \times \RR^N$ (this will always be
clear from the context), with the flat measure. When $\Theta >0$ is a
given positive locally locally integrable weight function, we write
$L^p(\Theta)$ for the Lebesgue space with measure $\Theta \, dx$, that
is the Banach space defined by the norm
$$
\| f \|_{L^p(\Theta)} := \left( \int_{\RR^N} |f|^p \, \Theta 
\right)^{1/p}.
$$
We shall sometimes use subscript letters in order to recall and
emphasize which variable we are considering for the Lebesgue space. 
\medskip

Our first result could thus read as follows:

\begin{theorem}\label{thm:0}
  Assume that $L$ is given by (\ref{eq:Lsimple}) and that $F$
  satisfies (\ref{eq:poly}) with $\alpha \in (0,2)$.  Assume
  furthermore that $f_0 \in L^2(F^{-1})$ 
  and let $f^\eps$ be the solution of (\ref{eq:eps1}) with
  $\theta(\eps)=\eps^\alpha$ and initial condition $f_0$.
 
  Then, when $\eps$ goes to zero, $f^\eps$ converges in
  $L^\infty(0,T;L^2(\RR^N\times\RR^N))$-weak to $\rho \, F$ with $\rho
  = \rho(t,x)$ the unique solution to the fractional diffusion
  equation (\ref{eq:D1}), (\ref{eq:D2}), (\ref{eq:D3}).
\end{theorem}

\smallskip

Note that under (\ref{eq:poly}), $\alpha$ can be characterized as
follows:
\begin{equation}\label{def:alphaASbdd}
 \alpha = \sup\left\{a \ge  0\, ;\, \int_{\RR^N} |v|^aF(v)\, dv <\infty \right\}.
\end{equation} 
This characterization will be more obvious in Section \ref{sectionMR}
when we add to $F$ a slowly varying function. In particular, we will
see that the fact that the moment of order $\alpha$ is finite or not
does not affect the asymptotic equation.

It is also worth pointing out that when $L$ is given by
(\ref{eq:Lsimple}), fractional diffusion is only observed when the
moment of order $2$ is unbounded, i.e. when the energy associated to
the equilibrium distribution function $F$ is infinite.  We will see,
however, in Theorem \ref{thm:1} that fractional diffusion may arise
even when the energy of $F$ (and higher order moments) is finite,
provided the collision frequency
$$\nu(v) = \int_{\RR^N} \sigma(v',v)\, dv'$$ 
in (\ref{eq:LL}) is degenerate for high velocities: $\nu(v)\sim
|v|^{\beta}$ for large $v$ with $\beta < 0$ small enough (for instance
with $\alpha > 2$ and $\beta < 2- \alpha <0$). 

Conversely, even when the energy of the equilibrium if infinite, it is
possible in some cases to recover a classical diffusion limit. This is
the case for instance when $\alpha \in (1,2)$ and $2-\alpha < \beta <
1$, or $\alpha \in (0,1)$ and $1 < \beta < 2 - \alpha$ (see
Theorem~\ref{thm:3}). These ranges of the parameters $\beta$ are quite
surprising as well as the fact that $\gamma(\alpha,\beta)$ (the order
of the limit fractional diffusion operator) is decreasing (to $1$) as
a function of $\beta$ when $\alpha \in (0,1)$ and increasing (to $1$)
as a function of $\beta$ when $\alpha \in (1,2)$, whereas one would
expect at first guess from physics that the stronger $\beta$, the
stronger the diffusion at the limit as it enhances collisions for high
velocities.  This calls for a satisfying physical interpretation.


\medskip

Theorem \ref{thm:0} states that in the long time limit, the particles
evolve according to an anomalous diffusion process (one calls a
diffusion process ÒanomalousÓ if the mean square displacement grows
like a nonlinear function of the time variable in the long time
limit).  Anomalous diffusion limits for kinetic models are well known
in the case of a gas confined between two plates, when the distance
between the plates goes to $0$ (see \cite{BGTh}, \cite{G}, \cite{D1},
\cite{D2}).  In that case, the limiting equation is still a standard
diffusion equation, but the time scale is anomalous ($\theta(\eps)\sim
\eps^2 \ln(\eps^{-1})$).  The particles travelling in directions
nearly parallel to the plates are responsible for the anomalous
scaling.  We will see that a similar behavior arises here when
$\alpha=2$ in (\ref{eq:poly}) (i.e.  when the second moment is
unbounded, but all moments of smaller orders are bounded), see Theorem
\ref{thm:3}.

\smallskip

A fractional diffusion equation has been obtained as a diffusive limit
from a linear phonon-Boltzmann equation simultaneously and
independently by M. Jara, T.  Komorowski and S. Olla in \cite{JKO}, by
a different probability approach.  Let us also mention a work in
progress \cite{MT} by B. Texier and the first author, where a
derivation of fractional diffusion equations from kinetic models is
obtained in the framework of a gas confined between two plates with
singular equilibrium distribution functions, by analytic means.

\medskip

Our approach relies on the use of Laplace-Fourier transforms and a
careful computation of the asymptotic behavior of the symbol of the
differential operator. This approach is quite simple and allows for an
explicit computation of the coefficients of the asymptotic equation
under optimal assumptions. Preliminary computations seem to indicate
that it can generalize to equations involving an external force field,
but it seems harder to generalize to a non-linear collision kernel
(Pauli's statistic for instance). A more ``non-linear'' approach, as
was the moments method in the classical diffusion case, would be
welcomed.

\smallskip

Let us now briefly outline the contents of the paper: The proof of
Theorem~\ref{thm:0}, which corresponds to the simplest case, is
presented in the next section. In Section~3, we state our main result
(Theorem~\ref{thm:1}) which generalizes Theorem~\ref{thm:0} to a large
class of linear collision operator $L$ of the form (\ref{eq:LL}). We
also address the critical case ($\alpha=2$) (Theorem~\ref{thm:2}) and
the classical diffusion case in Theorem~\ref{thm:3}.
Theorem~\ref{thm:1} (as well as Theorem~\ref{thm:3}) are then proven
in Section \ref{sec:4} while Section \ref{sec:5} is devoted to the
proof of Theorem~\ref{thm:2}.


\section{Proof of Theorem~\ref{thm:0} (simplest scenario)}\label{sec:mainidea}


The goal of the section is to present the main ideas of the proof,
which are in fact quite simple.  We thus make all sort of simplifying
assumptions in order to focus on the important aspect of the proofs.
We recall that $L$ is given by (\ref{eq:Lsimple}), and we assume,
instead of (\ref{eq:poly}), that the equilibrium function $F$
satisfies:
\begin{equation}\label{hypFsimple}
F(v) \le \kappa_0 \, |v|^{-N-\alpha} \mbox{ for all $v \in \RR^N$}, 
\quad
F(v) = \kappa_0 \, |v|^{-N-\alpha}\quad\hbox{if}\quad |v| \ge 1,
\end{equation} 
with $\kappa_0 > 0$. 

The proof is divided in 4 steps. Note that the proof of the main
result in Section \ref{sec:4} will follow exactly the same steps.

\medskip\noindent
{\bf Step 1: A priori estimates for $f^\eps$. } 

\noindent
The solution $f^\eps$ of (\ref{eq:eps1}) (with
$\theta(\eps)=\eps^\alpha$) satisfies
\begin{eqnarray*}
\eps^\alpha \, {d \over dt} \int_{\RR^{2N}} {(f^\eps)^2 \over 2} \, F^{-1} \, dvdx 
&=&  \int_{\RR^{2N}} L(f^\eps) \, f^\eps \, F^{-1}  \, dvdx \\
&=&  \int_{\RR^{2N}} [ (\rho^\eps)^2 \, F - (f^\eps)^2 \, F^{-1} ]  \, dvdx \\
&=& -  \int_{\RR^{2N}} [ f^\eps -  \rho^\eps \, F]^2 \, F^{-1}  \, dvdx, 
\end{eqnarray*}
from which we deduce the two estimates
\begin{equation}\label{bdd:thm01}
\sup_{t \ge 0} \int_{\RR^{2N}} {(f^\eps(t,.))^2 \over F}  \, dvdx \le 
 \int_{\RR^{2N}}  {f_0^2 \over F}  \, dv \, dx =  \| f_0 \|_{L^2(F^{-1})},
\end{equation}
and 
\begin{equation}\label{bdd:thm02}
\int_0^\infty \int_{\RR^{2N}} [ f^\eps -  \rho^\eps \, F]^2 \, F^{-1}
\, dv \, dx \, dt \le   
\frac{\eps^\alpha}2  \,  \| f_0 \|_{L^2(F^{-1})}.
\end{equation}
Cauchy-Schwarz inequality also gives:
$$
\rho^\eps (t,x) = \int_{\RR^N} {f^\eps \over F^{1/2}}Ê\, F^{1/2} \, dv 
\le \left( \int_{\RR^N} {(f^\eps)^2 \over F}Ê\, dv \right)^{1/2},
$$
so that $\rho^\eps(t,x)$, as well as $L(f^\eps)$, are
well defined a.e., and
\begin{equation}\label{bdd:thm03}
\sup_{t \ge 0} \int_{\RR^{N}}  \rho^\eps (t,.)^2 \, dx \le 
  \| f_0 \|_{L^2(F^{-1})}. 
\end{equation}

\smallskip\noindent {\bf Step 2: Another formulation of the rescaled
  equation. }

\noindent
We denote by $\widehat {f^\eps}$ the Laplace-Fourier transform of
$f^\eps$ with respect to $t$ and $x$, defined by
$$
 \widehat {f^\eps} (p,k,v) = \int_{\RR^N} \int_{0}^{\infty} e^{-pt} \,
 e ^{-i k\cdot x } \, f^\eps(t,x,v)\, dt\, dx,
 \quad p > 0, \,\, k \in \RR^N.$$ 
The function $\widehat {f^\eps}$ then satisfies
$$  \eps^\alpha  \, p \, \widehat {f^\eps} - \eps^\alpha \, \widehat {
  f_0} +\eps  \, i \, v\cdot k \widehat {f^\eps} 
=   \langle  \widehat {f^\eps} \rangle \, F- \widehat {f^\eps} ,
$$
where $\widehat {f_0}$ denotes the Fourier transform of the initial
datum $f_0$. We can rewrite this equality as
$$
\widehat {f^\eps}  =  \frac{ F  }{1+  \eps^{\alpha} \, p + \eps \, i
  \, v\cdot k  } \,  \widehat{\rho^\eps}
+  \frac{ \eps^{\alpha} \, \widehat{f_0}}{ 1+  \eps^{\alpha} \, p +
  \eps \, i \, v\cdot k   }\, ,
$$
with $\widehat {\rho^\eps}(p,k)  = \langle  \widehat{f^\eps} \rangle(p,k)$ the
Laplace-Fourier transform of $\rho^\eps$, and 
integrating this equality with respect to $v$, we obtain:
$$
\widehat {\rho^\eps}  = \left( \int_{\RR^N}  \frac{ F(v)}{1+
    \eps^{\alpha} \, p + \eps \, i \, v\cdot k  } \, dv \right) \, \widehat {\rho^\eps} 
+ \left( \int_{\RR^N} \frac{ \eps^{\alpha} \, \widehat{f_0}  }{ 1+
    \eps^{\alpha} \, p + \eps \, i \, v\cdot k  }\, dv \right).
$$
The normalization condition for $F$ now yields:
\begin{equation}\label{eq:10}
  \int_{\RR^N} \frac{ \widehat{f_0}}{ 1+ \eps^{\alpha} \, p 
     + \eps \, i \, v\cdot k   } \, dv+ a^\eps \, \widehat {\rho^\eps} 
  = 0, 
\end{equation}
with
$$ 
a^\eps(p,k):= \frac{1}{ \eps^{\alpha}}\int_{\RR^N}  
            \!\! \left( \frac{ 1}{1+  \eps^{\alpha} \,p + \eps \, i \,
                v\cdot k  } -1\right) \! F(v) \, dv.
$$
Let us prove that that the first term converges to $\widehat{\rho_0}$
when $\eps$ goes to zero.  The assumption $f_0 \in L^2(F^{-1})$
implies in particular that $f_0 \in L^2_x (L^1 _v)$. Hence its the
Fourier transform $\widehat f_0$ also belongs to $L^2_k (L^1 _v)$ by
Parseval equality, which means that $\widehat f_0$ is integrable in
$v$ for almost all $k$.  This allows to apply the Lebesgue dominated
convergence theorem, which yields, for almost every $k$,
$$
\int_{\RR^N}   \frac{ \widehat{f_0}}{ 1+ \eps^{\alpha} \, p 
     + \eps \, i \, v\cdot k   }  \, dv 
    \longrightarrow \int_{\RR^N}  \widehat f_0 \, dv = \widehat \rho_0.
$$

So we are left with the task of studying the limit, as $\eps$
goes to zero, of the coefficient $a^\eps$.

\medskip\noindent {\bf Step 3: The cornerstone argument of the proof:
  where fractional diffusion symbol appears. }

\noindent
A simple computation leads to
\begin{eqnarray*}
a^\eps(p,k )
& = &   
-  \frac{1}{ \eps^{\alpha}}\int_{\RR^N} 
  \frac{  \eps^{\alpha} p + \eps iv\cdot k}{1+  \eps^{\alpha} p + \eps iv\cdot k   } \, F(v) \, dv  \\
&=& - p \int_{\RR^N} 
   \frac{ 1+\eps^{\alpha} p}{ (1+  \eps^{\alpha} p)^2 + \eps^2 (v\cdot k  )^2 } \,  F(v) \, dv \\
&& 
-\frac{1}{\eps^\alpha}\int_{\RR^N} 
   \frac{ (\eps v\cdot k)^2}{ (1+  \eps^{\alpha} p)^2 + \eps^2 (v\cdot k  )^2 } \,  F(v) \, dv
\end{eqnarray*}
(the term involving $\eps i v \cdot k$ on the numerator vanishes
thanks to the symmetry of $F$). The first term in the right hand side
is bounded by $|p|$ (uniformly in $\eps$) and the dominated
convergence theorem readily implies that it converges to $-p
\int_{\RR^N} F(v)\, dv = -p$ as $\eps$ goes to zero. So it only
remains to study
$$
d^\eps(p,k) := \int_{\RR^N} 
  \frac{ \eps^{2-\alpha} (v\cdot k)^2}{ (1+  \eps^{\alpha} p)^2 + \eps^2 (v\cdot k  )^2 }  F(v) \, dv .
$$

If the second moment of $F$ is bounded, then it is readily seen that
with $\alpha=2$, we have
$$ 
d^\eps(p,k) \longrightarrow \int_{\RR^N} (v\cdot k)^2\, F(v)\, dv 
     =  \kappa \, |k|^2 ,\quad  \kappa \in (0,\infty),
$$
and the limit $\eps\rightarrow 0$ in (\ref{eq:0}) leads to the
diffusion equation (\ref{eq:diffusion}).

\smallskip

With (\ref{hypFsimple}) and $\alpha \in (0,2)$, the second moment of
$F$ is unbounded.  We then claim that for any $p \ge 0$, $k \in
\RR^N$, we have:
\begin{equation}\label{cvgce:deps}
| d^\eps(p,k) | \le \kappa \, |k|^\alpha 
\quad \hbox{and} \quad 
d^\eps(p,k) \mathop{\longrightarrow}_{\eps \to 0}    \kappa \, |k|^\alpha ,
\end{equation}
with $  \kappa \in (0,\infty)$ given by (\ref{eq:D3}). 

As a matter of fact, the first inequality in (\ref{cvgce:deps}) follows from (\ref{hypFsimple}):
$$
0 \le d^\eps(p,k) \le  \int_{\RR^N} \frac{ \eps^{2-\alpha} (v\cdot k)^2}{ 1 + \eps^2 (v\cdot k  )^2 } 
\, {\kappa_0 \over  |v|^{N+\alpha}} \, dv = \kappa \, |k|^\alpha,
$$
where the last equality is obtained
by making the change of variables $w:= \eps \, |k| \, v$. 
And in order to get the convergence of $d^\eps$, we simply write 
$$
d^\eps(p,k) = d^\eps_1(p,k) + d^\eps_2(p,k) 
$$
with
\begin{eqnarray*}
d^\eps_1(p,k)  & = &  \int_{|v| \le 1}  \ \frac{ \eps^{2-\alpha} (v\cdot k)^2}{ (1+  \eps^{\alpha} p)^2 + \eps^2 (v\cdot k  )^2 }  F(v) \, dv \\
 & \le  & \eps^{2-\alpha} \int_{|v| \le 1}   (v \cdot k)^2 \,  F(v) \, dv \\
 & \le  & \eps^{2-\alpha}|k|^2 \int_{|v| \le 1}    F(v) \, dv \to 0
\end{eqnarray*}
and
\begin{eqnarray*}
d^\eps_2(p,k)  & = &  \int_{|v| \ge 1}  \ \frac{ \eps^{2-\alpha} (v\cdot k)^2}{ (1+  \eps^{\alpha} p)^2 + \eps^2 (v\cdot k  )^2 }  \, {\kappa_0 \over  |v|^{N+\alpha}} \,dv \\
& = & |k|^\alpha  \int_{|w| \ge \eps \, |k| }  \ \frac{ w_1^2}{ (1+  \eps^{\alpha} p)^2 + w_1^2 }  \, {\kappa_0 \over  |w|^{N+\alpha}} \,dw 
\to \kappa \, |k|^\alpha,
\end{eqnarray*}
where we make again the change of variables $w:= \eps \, |k| \, v$ and  use the dominated convergence theorem. 
We have thus shown:
\begin{proposition}
If $\alpha\in(0,2)$, then
$$ a^\eps(p,k) \longrightarrow -p-\kappa |k|^\alpha \quad \mbox{ as } \eps\to 0$$
with  $  \kappa \in (0,\infty)$ given by (\ref{eq:D3}).
Furthermore, $ a^\eps(p,k)$ satisfies
$$ |a^\eps(p,k)| \leq |p|+\kappa|k|^\alpha.$$
\end{proposition}

\medskip\noindent
{\bf Step 4: Conclusion. } 

\noindent 
From the bound (\ref{bdd:thm03}) and up to extraction of a
subsequence, we know that there exists $\eta \in
L^\infty(0,\infty;L^2(\RR^N))$ such that $\rho^\eps \wto \eta$ weakly
in $L^\infty(0,\infty;$ $L^2(\RR^N))$.  On the one hand, from
(\ref{bdd:thm03}) again, we have $ \widehat{\rho^\eps}$ bounded in
$L^\infty(a,\infty;L^2(\RR^N))$ for any $a > 0$. On the other hand,
for any $\varphi = \varphi(t) \in \DD(0,\infty)$ its Laplace transform
$\LL \, \varphi$ belongs to $L^1(0,\infty)$ and for any $\psi =
\psi(x) \in \SS(\RR^N)$ its Fourier transform $\FF \, \varphi$ belongs
to $\SS(\RR^N)$ so that
$$
\langle  \widehat{\rho^\eps}, \varphi \otimes \psi \rangle = 
\langle \rho^\eps,\LL\varphi \otimes \FF\psi \rangle \to 
\langle  \eta, \LL\varphi \otimes \FF\psi \rangle 
= \langle  \widehat{\eta}, \varphi \otimes \psi \rangle
$$
as $\eps \to 0$.  We easily deduce that $ \widehat{\rho^\eps} \wto
\widehat{\eta}$ weakly in $L^\infty(a,\infty;L^2(\RR^N))$ for any $a >
0$.  Gathering all the convergence results established above, we may
pass to the limit in (\ref{eq:10}) and we get that $\eta$ satisfies
$$
 \widehat{\rho_0} +(-p- \kappa|k|^\alpha) \, \widehat \eta 
 = 0 \quad \hbox{for a.e.} \,\, p > 0, \,\, k \in \RR^N.
$$
Finally the unique solution $\rho$ to the fractional diffusion
equation (\ref{eq:D1}), (\ref{eq:D2}), (\ref{eq:D3}) also satisfies
 \begin{equation}\label{eq:laplace-fourier}
 \widehat{\rho_0} +(-p- \kappa|k|^\alpha) \, \widehat \rho 
 = 0 \quad \hbox{for a.e.} \,\, p > 0, \,\, k \in \RR^N,
 \end{equation}
so that $ \widehat \eta = \widehat \rho\,$ a.e., and then $\eta =
\rho$ because the Laplace-Fourier transform is a one-to-one
mapping (say in ${\cal S}' ([0,\infty) \times \RR^N)$).
We conclude that $f^\eps \wto \rho \, F$ weakly in $L^\infty(0,T;
L^2(\RR^N \times \RR^N))$ thanks to (\ref{bdd:thm02}). \qed

\vspace{10pt}


\section{Diffusion limit  results for general collision operator}\label{sectionMR}


Having exposed the main ideas of the proof in the previous section, we
now come back to more general collision operators, and state our main
result. 

We recall that $L$ now is of the form
\begin{eqnarray*}
L(f)  & = &  \int_{\RR^N} \Big[ \sigma(v,v') f(v') - \sigma(v',v) f(v)\Big]\, dv'\\[3pt]
 & = & K(f) - \nu f
\end{eqnarray*}
with 
$$ 
K(f) =   \int_{\RR^N} \sigma(v,v') f(v') \, dv' , \qquad \nu(v) =  \int_{\RR^N} \sigma(v',v) \, dv' .
$$
Such an operator is obviously linear, well-defined as a (possibly)
unbounded operator with domain $L^1(\nu)$, and closed (in fact it is
bounded from $L^1(\nu)$ to $L^1)$. It is also straightforwardly conservative:
$$ \int_{\RR^N} L(h)\, dv = 0 \quad \mbox{ for all } h \in L^1(\nu).$$
The choice of the cross-section $\sigma(v,v')$ is crucial. We start
with the following structural assumption.

\medskip
\noindent {\bf Assumptions (A1)} {\it The cross-section $\sigma$
  is locally integrable on $\RR^{2N}$, non negative and the collision
  frequency $\nu$ is locally integrable on $\RR^N$ and satisfies
$$ 
\nu(-v)=\nu(v) > 0 \quad \mbox{ for all } v\in\RR^N.
$$

\noindent {\bf Assumptions (A2)} There exists a function $0 \le F \in
L^1(\nu)$ such that $|v|^2 \, \nu(v)^{-1} \, F$ is locally integrable
and 
\begin{eqnarray}\label{def:F}
\nu (v) \, F(v) = K(F)(v) =   \int_{\RR^N} \sigma(v,v') F(v') \, dv',
\end{eqnarray}
which means that  $F$ is a equilibrium distribution (i.e.  $L(F)=0$). 
Furthermore, the function $F$ is symmetric,  positive and normalized to $1$: 
$$ 
F(-v) =F(v) > 0 \,\, \mbox{ for all }  \,\, v\in\RR^N
 \quad\hbox{and}\quad \int_{\RR^N } F(v)\, dv=1.
$$
}

\medskip

Note that under classical assumptions on $\sigma$, the existence of an
equilibrium function is in fact a consequence of Krein-Rutman's
theorem (see \cite{DGP} for details).  A particular case in which this
condition is satisfied is when $\sigma$ is such that 
\begin{eqnarray}\label{def:micro}
\forall \, v,v' \in \RR^N \qquad 
\sigma(v,v') =b(v,v') F(v), \quad  b(v',v) = b(v,v'),
\end{eqnarray} 
for some $b \in L^1 _{loc}(\RR^{2N})$. In that case, we say that
$\sigma$ satisfies a {\it detailed balanced principle} or a {\it
  micro-reversibility principle}, while the more general assumption
(\ref{def:F}) is called a {\it general balanced principle}.
  
\medskip

Next, we need to make precise the behavior of $F$ and $\nu$ for large
$|v|$. For that purpose, we recall that a slowly varying function is a
measurable function $\ell:\RR_+\rightarrow \RR$ such that
$$ \ell(\lambda s) \sim \ell (s) \quad \mbox{ as } s\rightarrow \infty \mbox{ for all } \lambda>0.$$
Example of slowly varying functions are positive constants, functions
that converge to positive constants, logarithms and iterated
logarithms.

In our main result, we assume that $F$ is a regularly varying function
of index $-(N+\alpha)$ with $\alpha>0$ and that $\nu$ behaves in the
large velocity asymptotic like a power function. More precisely, we
make the following assumption:


\medskip
\noindent {\bf Assumptions (B1)}
{\it 
 There exists $\alpha >0$ and a slowly varying function $\ell$ such that
\begin{equation}\label{eq:Fdef}
F(v) = F_0(v) \ell(|v|),
\end{equation}
where $F_0$ is such that
\begin{equation}\label{eq:F0def}
|v|^{\alpha+N} F_0(v) \longrightarrow  \kappa_0 \in (0,\infty) \  \mbox{ as } |v|\rightarrow \infty.
\end{equation}

\noindent
{\bf Assumptions (B2)} 
There exists $\beta\in\RR$ and a positive constant $\nu_0$ such that
\begin{equation}\label{eq:nu0def}
 |v|^{-\beta} \nu(v) \longrightarrow \nu_0 \ \mbox{ as } |v|\rightarrow \infty.
\end{equation}

\noindent
{\bf Assumptions (B3)}   Finally we assume that there exists a constant $M$  such that
\begin{equation}\label{bdd:bnuF}
\int_{\RR^N} F'  \, {\nu \over b} \, dv' 
+ \left( \int_{\RR^N} {F' \over \nu'} \, {b^{2} \over \nu^2} \, dv' \right)^{1/2} 
\le M \quad\mbox{ for all }\quad v \in \RR^N,
\end{equation}
with $b = b(v,v') := \sigma(v,v') \, F^{-1}(v)$. 
}

\begin{rk&ex}\label{rk1}
\item[(i)] The constant $ \kappa_0$ in (\ref{eq:F0def}) could actually
  be a function of the direction $v/|v|$ without any additional
  difficulties. We will take $ \kappa_0$ constant in order to keep
  things simple.
\item[(ii)] The condition (B2)  implies that the collision frequency satisfies:
\begin{equation}\label{eq:nub}
  \nu_1 \, \langle v\rangle^\beta \leq \nu(v) \leq \nu_2 \, \langle v
  \rangle^\beta \ \mbox{ for large } v\in\RR^N.
 \end{equation}
 for some constants $\nu_1,\nu_2 \in (0,\infty)$.
\item[(iii)] The conditions (B2) and (B3) are fulfilled for a
  collision kernel $\sigma$ satisfying the detailed balance principle
  (\ref{def:micro}) where $b$ satisfies
$$
b(v,v') = \langle v \rangle^\beta \, \langle v' \rangle^\beta  
\quad\hbox{or}\quad b(v,v') =   \langle v-v' \rangle^\beta, \quad \beta < \alpha,
$$
since then $\nu(v) \sim \langle v \rangle^\beta$ in both cases (for
the second example, we refer to the proof of Lemma~\ref{lem:LA} in the
appendix where the main arguments of the proof of that statement is
presented). Let us notice that such collision kernels satisfy
\cite[Assumptions (A1), (A2), (A3)]{DGP} but of course not in general
\cite[Assumptions (B2), (B3)]{DGP} except when $\beta > 2-\alpha$.
\item[(iv)] The conditions (B2) and (B3) are also fulfilled for a
  (more physical) collision kernel $\sigma$ satisfying the detailed
  balance principle (\ref{def:micro}) with
$$
b(v,v') = |v-v'|^\beta
$$
and under the additional restriction $\beta \in
(-\min\{\alpha;N/2\},\min\{\alpha;N)\}$. We refer to Lemma~\ref{lem:LA} in
the appendix where the proof of that statement is presented.
\item[(v)] Our assumption (B3) is a bit more general than the
  corresponding assumptions (A3) in \cite{DGP} since for instance it
  is fulfilled by the collision kernels of the point (iv) above with
  $\beta \in (0,\min\{\alpha;N\})$, while such a collision kernel does
  not satisfy \cite[Assumption (A3)]{DGP}.
\end{rk&ex}

We can now  state our main theorem:

\begin{theorem}[Fractional diffusion limit]\label{thm:1} 
  Assume that Assumptions (A1-A2) and (B1-B2-B3) hold with $\alpha>0$
  and $\beta <\min\{\alpha;2-\alpha\}$.  Define
$$
\gamma := \frac{\alpha-\beta}{1-\beta}, \quad \mbox{ and
}\quad\theta(\eps):= \ell(\eps^{-\frac{1}{1-\beta}}) \, \eps^\gamma.
$$
It is worth noticing that we have here $\beta < 1$ as well $\gamma <
2$ for these ranges of the parameters $\alpha$ and $\beta$.
Assume furthermore that $f_0 \in L^2(F^{-1})$ 
and let $f^\eps$ be the solution of (\ref{eq:eps1}), with that choice of 
$\theta$ and initial data $f_0$.

Then, $(f^\eps)$ converges in
$L^\infty(0,T;L^2(\RR^N\times\RR^N))$-weak to a function $\rho(t,x) \,
F(v)$ where $\rho(t,x) $ is the unique solution of the fractional
diffusion equation of order $\gamma$:
\begin{eqnarray}
  && \pa_t \rho + \kappa \, (-\Delta_{x})^{\gamma/2} \rho= 0 
 \qquad \mbox{ in }(0,\infty) \times \RR^N,\label{eq:D12} \\[4pt]
  && \rho(0,.) = \rho_0  \quad \qquad  \qquad\,  \qquad \mbox{ in }\RR^N,\label{eq:D22}
\end{eqnarray}
with $\kappa$ given by
\begin{equation}\label{eq:D4}
  \kappa = \frac{ \kappa_0 \, \nu_0}{1-\beta} 
     \int_{\RR^N} \frac{ w_1^2}{ \nu_0^2 + w_1^2 } \, {1 \over  |w|^{N+\gamma}} \, dw. 
\end{equation}

\end{theorem}

\begin{remark}\label{rk2}
\item[(i)] Note that in Theorem~\ref{thm:1} we always have $\beta <1$
from the assumptions, and the condition $\beta
<\min\{\alpha;2-\alpha\}$ is equivalent to the condition $\gamma <2$.

\item[(ii)] When $\beta=0$ (i.e. $\nu$ bounded below by a positive
  constant), then we have $\gamma=\alpha$, and we recover a result
  similar to that of Theorem \ref{thm:0} with more general collision
  operator and equilibrium states. Note in particular that in that
  case, the addition of a slowly varying part in $F$ has modified the
  time scale $\theta(\eps)$, but not the limiting equation.  
  
\item[(iii)] When $\alpha \in (0,1)$, the function $\beta \in
  (-\infty,\alpha) \mapsto \gamma(\beta) = (\alpha - \beta) /
  (1-\beta)$ is decreasing from $1$ to $0$. When $\alpha =1$, it is
  defined on $(-\infty,1)$ and identically constant to $1$. When
  $\alpha >1$, the function $\beta \in (-\infty,2-\alpha) \mapsto
  \gamma(\beta) = (\alpha - \beta) / (1-\beta)$ is increasing from $1$
  to $2$. It is thus always possible to obtain fractional diffusion
  limit for a kinetic collisional equation with regularly varying
  equilibrium $F$ of index $-(N+\alpha)$. In particular, fractional
  diffusion behavior can arise even when $F$ has finite energy
  ($\alpha>2$).

\end{remark}


The case $\gamma=2$ (which may occur when $\alpha > 1$) is critical in
the sense that even though the second moment may be infinite (and the
usual method yield an infinite diffusion coefficient), the asymptotic
behavior can still be described by a standard diffusion equation under
the appropriate time scale.  In this case, the exact behavior of the
slowly varying function $\ell$ is crucial, as it may determine whether
the second moment is finite or infinite.  When it is finite, (for
instance $\ell(|v|) =(\ln|v|)^{-2}$), then the usual technic yields
the diffusion equation (\ref{eq:diffusion}) under the classical time
scale $\theta(\eps)=\eps^2$.  When it is infinite, then the asymptotic
behavior is still described by a diffusion equation, but the time
scale has to be modified. These two situations are
included in the next theorems: the first one corresponds to anomalous
diffusive time scales, while the second one corresponds to classical
diffusive time scales.

\begin{theorem}[Classical diffusion limit with anomalous time scale]
\label{thm:2} 
Assume that Assumptions (A1-A2) and (B1-B2-B3) hold with
$$
\alpha > 1 \quad \hbox{and} \quad
\beta = 2 - \alpha \quad (\mbox{i.e.} \ \gamma =2),
$$
and $\ell$ such that 
\begin{equation}\label{ellcritic}
\ell(r) \, \ln (r) \to + \infty \quad \mbox{ as } r \to + \infty
\end{equation}
(note that this implies in particular that the second moment of $F$ is
infinite). 

Then define 
$$
\theta(\eps) = \eps^2 \, \ell(\eps^{-\frac{1}{1-\beta}}) \, \ln(\eps^{-1}).
$$
Assume furthermore that $f_0 \in L^2(F^{-1})$  
and let $f^\eps$ be the solution of (\ref{eq:eps1}), with
$\theta(\eps)$ defined as above and initial data $f_0$.

Then, $(f^\eps)$ converges in
$L^\infty(0,T;L^2(\RR^N\times\RR^N))$-weak to $\rho \, F$ where $\rho
= \rho(t,x)$ is the unique solution to the standard diffusion equation
$$ \pa_t \rho - \kappa \, \Delta_x \rho =0$$
with $\kappa$ given by 
$$\kappa = \frac{\kappa_0 \, \nu_0}{(1-\beta)} \, 
\lim_{\lambda\to 0} \frac{1}{\ln(\lambda^{-1})}  
 \int_{|w|\geq \lambda} \frac{ w_1^2}{ \nu_0^2 + w_1^2 } \, {1 \over  |w|^{N+2}} \, dw. 
$$
\end{theorem}

\begin{remark}
\item[(i)] This is a critical case, in which the variance of the
  equilibrium distribution $F$ (and therefore the classical diffusion
  coefficient (\ref{eq:coefdiffclass})) are infinite, but the
  asymptotic symbol still is of order $2$. This case is still referred
  to as anomalous diffusion, even though we recover a standard
  diffusion equation, because the time scale is not the usual time
  scale.
\item[(ii)] Proceeding as Section \ref{sec:mainidea}, with $\alpha=2$,
  $\beta=0$, and $\ell \equiv 1$, we immediately see that if we would
  take a time scale $\theta(\eps)=\eps^2$, then we should expect to
  find for the diffusion coefficient:
$$
\kappa= \int_{\RR^N } \ \frac{ w_1^2}{ 1 + w_1^2 } \, {\kappa_0 \over
  |w|^{N+2}} \,dw = \infty.$$ This is why a different time scale has
to be considered in this limiting case $\alpha=2$. 
\end{remark}

\begin{theorem}[Classical diffusion limit with classical time scale] 
\label{thm:3} 
Assume that Assumptions (A1-A2) hold as well as the following bounds
\begin{equation}\label{bdd:g2}
  \int_{\RR^N} \left(  {\nu(v)  \over b(v,v')} + {|v'|^2 \over
      \nu(v')} \right) 
            \, F' \, dv' \le M \qquad \forall \, v \in \RR^N.
\end{equation}
Assume furthermore that $f_0 \in L^2(F^{-1})$ and let $f^\eps$ be the
solution of (\ref{eq:eps1}), with $\theta(\eps) = \eps^2$ and initial
data $f_0$.

Then $(f^\eps)$ converges in
$L^\infty(0,T;L^2(\RR^N\times\RR^N))$-weak and in $L^2((0,T) \times
\RR^N\times\RR^N))$-strong to a function $\rho \, F$ where $\rho =
\rho(t,x)$ is the unique solution of the standard diffusion equation
(\ref{eq:diffusion}) explained in the introduction (with the same
constant).
\end{theorem}

\begin{remark}\label{rk3}
\item[(i)] Remark that Asumptions (B1-B2-B3) with any $\beta >
  2-\alpha$ (i.e. $\gamma>2$) obviously implies (\ref{bdd:g2}), so
  that all the values of $\alpha > 0$, $\beta <\min\{1;\alpha\}$ are
  addressed by Theorems~\ref{thm:1}, \ref{thm:2} and \ref{thm:3}. One
  also sees from these exemples of coefficients that this theorem
  covers cases where the velocity distribution has infinite variance.
  \item[(ii)] Theorem~\ref{thm:3} is (in some direction!) slightly
    more general than \cite[Theorem 1]{DGP} since for instance
    \cite[Asumption (A3)]{DGP} implies that $b$ is uniformly bounded
    from below, whereas Theorem~\ref{thm:3} applies to the collision
    kernels of Remarks \& Examples~\ref{rk1} (iv) with $\beta \in
    (\max\{2-\alpha;0\},\min\{\alpha;N\})$, while they are not dealt
    with by \cite[Theorem 1]{DGP}.
\end{remark}

Figure 1  summarizes a large part of the results
presented in the last three theorems.

\begin{figure}\label{lafigure}
  \begin{center}
    \mbox{ \includegraphics{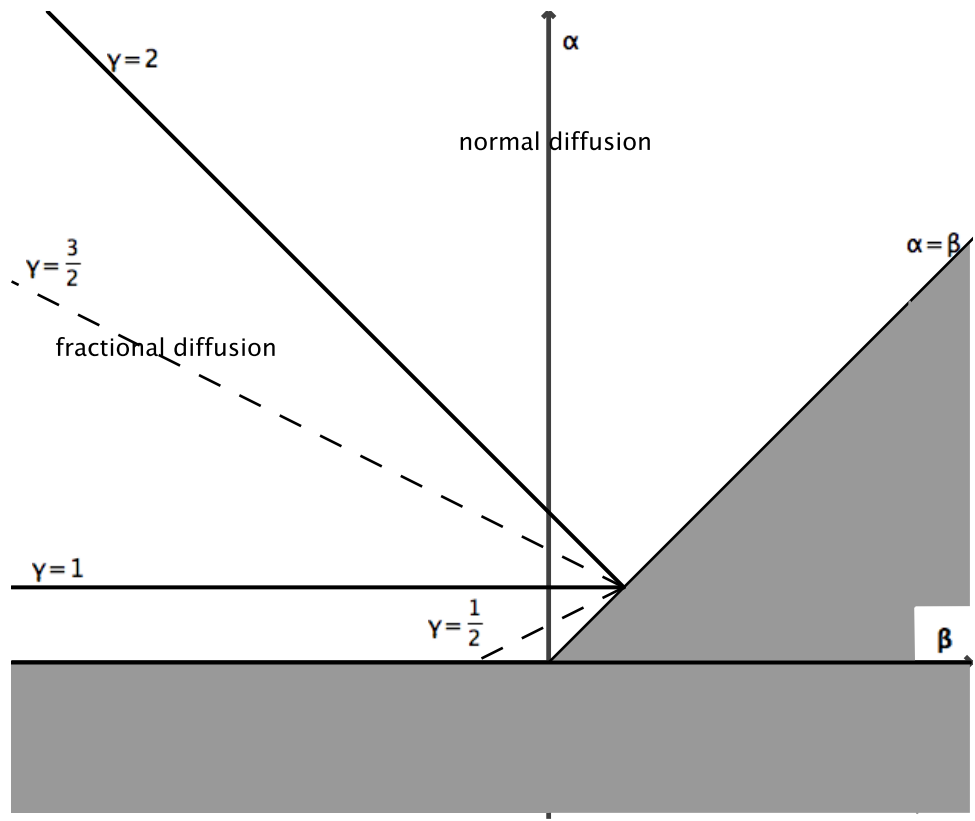} }
    \caption{Summary of the main results}
  \end{center}
\end{figure}

We now turn to the proofs of Theorems \ref{thm:1} and \ref{thm:2},
presented respectively in Sections \ref{sec:4} and \ref{sec:5}.
Concerning Theorem \ref{thm:3} we shall not write a full proof since
it just follows from the proof of \cite[Theorem 1]{DGP} by replacing
the proof of the auxiliary result in \cite[Proposition 1]{DGP} by
Lemma~\ref{lem:L} below.


\section{Proof of Theorem~\ref{thm:1}}\label{sec:4}


The general idea of the proof is the same as that of Theorem
\ref{thm:0} (see Section \ref{sec:mainidea}).  The main difference is
in the derivation of formulation (\ref{eq:10}): It is no longer
possible to work with the density $\rho^\eps$ directly, and we use
$K(f^\eps)$ instead leading to the corresponding formulation
(\ref{eq:pp}).  This leads to another difficulty, since (\ref{eq:pp})
involves an additional term which we have to show is of smaller order.
Finally, we will see that the computation of the asymptotic symbol is
a little bit more complicated than in Section \ref{sec:mainidea}
because of the collision frequency being velocity dependent, and
because of the presence of the slowly varying function $\ell$ in $F$
(this is the object of Proposition \ref{lem:aa}, in which we see
appearing the importance of the different time scale that takes into
account $\ell$ and $\beta$). 

Throughout this section (and the next one), we denote
$$\vphi(\eps): =\ell(\eps^{-\frac{1}{1-\beta}} ).$$

\subsection{A priori estimates}


\vspace{10pt}

The next lemma summarizes the key properties of the collision  operator $L$:
\begin{lemma}\label{lem:L}
  The operator $\frac{1}{\nu} L$ is bounded in $L^2(\nu F^{-1})$ and
  satisfies:
\begin{equation}\label{eq:Q}
  \int _{\RR^N} L(f) \frac{f}{F}\, dv \leq -\frac{1}{2M} 
         \int_{\RR^N} |f-\langle f\rangle F|^2\frac{\nu}{F}\, dv
                      \quad \mbox{ for all } f\in L^2(\nu F^{-1})
\end{equation}
where $\langle f\rangle = \int_{\RR^N} f(v)\, dv$ and $M$ is defined
in (\ref{bdd:bnuF}).
\end{lemma}

\medskip\noindent{\bf Proof of Lemma~\ref{lem:L}.} We adapt the proof
of \cite[Proposition 1 \& 2]{DGP}.

To show that $\frac{1}{\nu}L$ is bounded in $L^2(\nu F^{-1})$, we
obviously only have to check that $\frac{1}{\nu} K$ is bounded.  Using
Cauchy-Schwartz inequality and the fact that $K(F)=\nu F$, we get:
\begin{eqnarray*}
\left\| \frac{1}{\nu}K(f)\right\|^2_{\LnF} & = &  
\int \frac{K(f)^2}{\nu F}\, dv \\
& \leq & \int \frac{1}{\nu (v) F(v)} \int \sigma(v,v') F(v') dv' 
          \int \sigma(v,v') \frac{f(v')^2}{F(v')}\, dv'\, dv \\
& \leq & \int \frac{1}{\nu(v) F(v)} \nu(v)F(v) 
           \int \sigma(v,v') \frac{f(v')^2}{F(v')}\, dv'\, dv \\
& \leq & \int \int \sigma(v,v') \frac{f(v')^2}{F(v')}\, dv'\, dv \\
& \leq &  \int \nu(v') \frac{f(v')^2}{F(v')}\, dv' = \|f\|^2_{L^2(\nu F^{-1})}.\\
\end{eqnarray*}

In order to prove (\ref{eq:Q}), we write
\begin{eqnarray*}
\int _{\RR^N} L(f) \frac{f}{F}\, dv  & =  & \int _{\RR^N}\int _{\RR^N} \sigma(v,v')f'\frac{f}{F}\, dv\, dv' -\int _{\RR^N} \nu(v)\frac{f^2}{F}\, dv\\
& =  & \int _{\RR^N}\int _{\RR^N} \sigma(v,v')F'\frac{f'}{F'}\frac{f}{F}\, dv\, dv' -\int _{\RR^N} \nu(v)\frac{f^2}{F}\, dv.
\end{eqnarray*}
Next, we note that the second term in the right hand side can be rewritten
\begin{eqnarray*}
\int _{\RR^N} \nu(v)\frac{f^2}{F}\, dv & = & \int _{\RR^N} \int_{\RR^N} \sigma(v',v) F\frac{f^2}{F^2}\, dv\, dv'\\
& = & \int _{\RR^N} \int_{\RR^N} \sigma(v,v') F'\frac{f'^2}{F'^2}\, dv\, dv,
\end{eqnarray*}
as well as  (using the fact that $\nu F=K(F)$)
\begin{eqnarray*}
\int _{\RR^N} \nu(v)\frac{f^2}{F}\, dv & = &\int _{\RR^N} K(F)\frac{f^2}{F^2}\, dv \\
& = & \int _{\RR^N} \int_{\RR^N} \sigma(v,v') F'\frac{f^2}{F^2}\, dv\, dv.
\end{eqnarray*}

We deduce
\begin{eqnarray}
\int _{\RR^N} L(f) \frac{f}{F}\, dv
& =  & -\frac{1}{2}\int _{\RR^N}\int _{\RR^N} \sigma(v,v')F'  \left[ \frac{f'}{F'}-\frac{f}{F}\right]^2\, dv\, dv' .
\label{coercif1}
\end{eqnarray}
Integrating (in the $v'$ variable) the identity
$$
f \, F' - f' \, F = \left(\frac{f}{F} - \frac{f'}{F'}\right) \,  F \, F'
$$
we get
$$
g = \int_{\RR^N} \left(\frac{f}{F} - \frac{f'}{F'}\right) \,  F \, F' \, dv'
$$
where $g=f-\langle f\rangle F$.  
The Cauchy-Schwarz inequality implies
$$
g^2 \le \left( 
\int_{\RR^N} \left(\frac{f}{F} - \frac{f'}{F'}\right)^2 \, \sigma \,
F'  \, dv' \right) \  
\left( \int_{\RR^N} { F^2 \over \sigma} \, F'  \, dv' \right),
$$
so that 
\begin{equation}\label{coercif2} 
\int_{\RR^N}Ê{g^2 \over F} \, \nu \, dv \le
 \left( \sup_{v \in \RR^N}\nu \int_{\RR^N} { F \over \sigma} \, F'  \,
   dv' \right) \, \left(
  \int_{\RR^N} \, \int_{\RR^N} \left(\frac{f}{F} -
    \frac{f'}{F'}\right)^2 \, \sigma \, F'  \, dv' \, dv \right).
\end{equation}
Gathering (\ref{bdd:bnuF}), (\ref{coercif1}) and (\ref{coercif2}) we obtain (\ref{eq:Q}). 
\qed

\medskip

Using Lemma \ref{lem:L}, we can prove the following estimate on
$f^\eps$:
\begin{lemma}\label{lem:bound}
  The solution $f^\eps$ of (\ref{eq:eps1}) is bounded in
  $L^\infty(0,\infty;L^2 (F^{-1}))$ uniformly with respect
  to $\eps$.  Furthermore, it satisfies:
  $$ f^\eps = \rho^\eps F(v) + g^\eps, $$
  where the density $ \rho^\eps = \int_{\RR^N} f^\eps\, dv$ and the
  function $g^\eps$ are such that
\begin{equation}\label{eq:rho} 
\|\rho^\eps\|_{L^\infty(0,\infty,L^2 )} \leq \| f_0\|_{L^2 (F^{-1})} 
\end{equation}
and
\begin{equation}\label{eq:gg} 
\|g^\eps\|_{L^2(0,\infty;L^2 (\nu F^{-1}))} \leq C \, \|f_0\|_{L^2 (F^{-1})} \, \theta(\eps)^{1/2}.
\end{equation}

\item In particular  $\rho^\eps$ converges $L^\infty(0,T;L^2)$-weak to $\rho$, and 
$f^\eps$ converges  $L^\infty(0,T;L^2 (F^{-1}))$-weak to  $f=\rho(t,x)F(v)$.
\end{lemma}

\medskip\noindent{\bf Proof of Lemma~\ref{lem:bound}.} 
Multiplying (\ref{eq:eps1}) by $f^\eps/F$, we get:
\begin{eqnarray*}
\frac{1}{2} \, \frac{d}{dt} \int_{\RR^{2N}}  |f^\eps|^2\frac{1}{F}\,
dx\,dv & = &  
\frac{1}{\theta(\eps) }  \, \int_{\RR^{2N}}  L(f^\eps) \frac{f^\eps}{F}\\
& \leq &  - \frac{1}{2M \theta(\eps)} \, \int_{\RR^{2N}}|f^\eps-\rho^\eps  F|^2 \frac{\nu}{F}\, dx\,dv 
\end{eqnarray*}
which gives:
\begin{eqnarray*}
&&  \frac12 \, \int_{\RR^{2N}}
|f^\eps(t,x,v)|^2\frac{1}{F}\, dx\,dv 
  +  \frac{1}{2 M \, \theta(\eps)} \, 
  \int_0^t \int_{\RR^{2N}}|f^\eps-\rho^\eps  F|^2 \frac{\nu}{F}\, dx\,dv\, ds \\
&& \qquad\qquad\qquad\qquad\qquad\qquad\qquad\qquad  \leq
\frac12 \, \int_{\RR^{2N}} |f_0(x,v)|^2\frac{1}{F}\, dx\,dv.
\end{eqnarray*}
This inequality shows that $f^\eps$ is bounded in
$L^\infty(0,\infty,L^2(F^{-1}))$. Furthermore, denoting $g^\eps=
f^\eps-\rho^\eps F$ we also get
$$ 
\int_0^t \int_{\RR^{2N}} |g^\eps|^2 \frac{\nu}{F}\, dx\,dv\, ds \leq 
C \,\|f_0\|_{L^2(F^{-1})} \, \theta(\eps)  .
$$
\medskip
 
Finally, Cauchy-Schwarz inequality implies:
\begin{eqnarray*}
\int_{\RR^N} |\rho^\eps|^2\, dx & = & \int_{\RR^N} \left|\int_{\RR^N} f^\eps\, dv\right|^2 \, dx \\
& \leq & \int_{\RR^{2N}}  |f^\eps|^2\frac{1}{F}\, dv \int_{\RR^N} F\, dv \, dx =
 \int_{\RR^{2N}} |f^\eps|^2\frac{1}{F}\, dv \, dx .
\end{eqnarray*}
\qed

\vspace{20pt}

\subsection{Another formulation of the rescaled equation}


Proceeding as in Section \ref{sec:mainidea}, we denote by $\widehat
{f^\eps}(p,k,v)$ the Laplace-Fourier transform of $f^\eps(t,x,v)$ with
respect to $t$ and $x$, defined by
$$   \widehat {f^\eps} (p,k,v) = \int_{\RR^N} \int_{0}^{\infty} e^{-pt} e ^{-i k\cdot x } f^\eps(t,x,v)\, dt\, dx.$$ 
We define $\widehat {\rho^\eps} (p,k)$ and $\widehat {g^\eps} (p,k,v)$
similarly, and we denote by $\widehat { f_0}(k,v)$ the Fourier
transform of $ f_0$ with respect to $x$.
\medskip

Now, taking the Laplace-Fourier transform in (\ref{eq:eps1}), it is
readily seen that $\widehat {f^\eps}$ satisfies
$$  
\theta(\eps) \, p\,  \widehat {f^\eps} 
  -\theta(\eps) \, \widehat { f_0} + i \, v\cdot k \, \widehat{f^\eps} 
         =   K(  \widehat {f^\eps}) - \nu \, \widehat {f^\eps},
$$
which easily yields
$$
\widehat {f^\eps}(v)  = \frac{  \theta(\eps)}{ \nu(v)+ \theta(\eps)   p + \eps iv\cdot k  }\widehat{f_0} +  \frac{ 1}{ \nu(v)+  \theta(\eps)   p + \eps iv\cdot k }  K(  \widehat {f^\eps}) .
$$

Multiplying this equality by $\sigma(w,v)$ and integrating with respect to $v$, we obtain
\begin{eqnarray*}
K(\widehat {f^\eps})(w) 
& = &  
\int_{\RR^N} \frac{\theta(\eps)  \sigma(w,v) }
   {\nu(v)+ \theta(\eps)   p + \eps iv\cdot k   }\, \widehat{f_0}(v)  \, dv \\
&& + \int_{\RR^N}  \frac{  \sigma(w,v)}{( \nu(v)+ \theta(\eps)  p + \eps iv\cdot k  ) } 
   \, K(  \widehat {f^\eps})(v)  \, dv.  
\end{eqnarray*}
Finally, integrating with respect to $w$, we get:
\begin{eqnarray*}
\int_{\RR^N}K(\widehat {f^\eps})(v)\, dv 
& = &  
\int_{\RR^N} \frac{\theta(\eps)   \nu(v)}{\nu(v)+\theta(\eps)   p + \eps iv\cdot k   }\,  \widehat{f_0} (v) \, dv \\
&& + \int_{\RR^N}  \frac{\nu(v)}{ \nu(v)+  \theta(\eps)
  p + \eps iv\cdot k  } \, K(  \widehat {f^\eps})(v)  \, dv,   
\end{eqnarray*}
and thus 
\begin{eqnarray*}
0 & = &  
\int_{\RR^N} \frac{\nu(v) }{\nu(v)+  \theta(\eps)   p + \eps iv\cdot k   }\,  \widehat{f_0}(v) \, dv \\
&&+ \frac{1}{\theta(\eps)} \, \left( \int_{\RR^N} \left[
    \frac{\nu(v)}{ \nu(v)+  \theta(\eps)   p + \eps
      iv\cdot k  } -1\right] \, K(  \widehat {f^\eps})(v)  \, dv \right). 
\end{eqnarray*}

\smallskip 

Next, using Lemma \ref{lem:bound}, we  write
$$ \widehat{ f^\eps} = \widehat {\rho^\eps} F  +\widehat {g^\eps}$$
which leads to (using the fact that $K(F)=\nu F$):
$$ K(  \widehat {f^\eps}) =  \widehat {\rho^\eps} \nu F + K(\widehat {g^\eps}).$$
We deduce:
\begin{eqnarray}
\!\!\!\!\!  0  & = &  
\int_{\RR^N} \frac{\nu(v)   }{\nu(v)+  \theta(\eps)   p + \eps iv\cdot
  k   }\, 
\widehat{f_0} \, dv \nonumber \\[5pt]
&&+ \frac{1}{\theta(\eps)} \, 
    \left( \int_{\RR^N} \left[ \frac{\nu(v)}{ \nu(v)+  \theta(\eps)  p + \eps iv\cdot k  } -1\right] \, \nu F \,
      dv\right) \, \widehat {\rho^\eps}  \nonumber\\[5pt]
&&+ \frac{1}{ \theta(\eps)} \, \left( \int_{\RR^N} 
      \left[ \frac{\nu(v)}{ \nu(v)+  \theta(\eps)   p +
          \eps iv\cdot k  } -1\right] \, K(  \widehat {g^\eps})  \, dv \right). \label{eq:pp}
\end{eqnarray}

The rest of the proof consists in passing to the limit
$\eps\rightarrow 0$ in (\ref{eq:pp}). 

In the next two subsections, we shall show that the last term vanishes
in Lemma~\ref{lem:aa} and we study in Lemma~\ref{lem:a} the limit of
the second term
$$ 
a^\eps(p,k):= \frac{1}{ \theta(\eps)} \, \left(
  \int_{\RR^N} \left[ \frac{\nu(v)}{ \nu(v)+  \theta(\eps) \, p 
+ \eps \, i \, v\cdot k  } -1\right] \, \nu \, F  \, dv \right)
$$
as $\eps$ goes to zero.  This limit will provide the Fourier-Laplace
symbol of the asymptotic equation, and is therefore the cornerstone of
the proof.

Let us prove that the first term converges to $\widehat \rho_0$
  \begin{equation}\label{eq:lF0}
    \int_{\RR^N}  \frac{ \nu(v)}{ \nu+  \theta(\eps)   \,  p +
      \eps \, i \, v\cdot k   }  \widehat f_0 \, dv 
    \longrightarrow \rho_0, \qquad \mbox{for almost all } k\in\RR^N,\, p\in\RR_+
  \end{equation}

  The assumption $f_0 \in L^2(F^{-1})$ implies in particular that $f_0
  \in L^2_x (L^1 _v)$. Hence its the Fourier transform $\widehat f_0$
  also belongs to $L^2_k (L^1 _v)$ by Parseval equality, which means
  that $\widehat f_0$ is integrable in $v$ for almost all $k$.
  Together with (\ref{bdd:nu-1}), this allows to apply the Lebesgue
  dominated convergence theorem, which yields, for almost every $k$,
$$
\int_{\RR^N}  \frac{ \nu(v)}{ \nu+  \theta(\eps)   \,  p +
      \eps \, i \, v\cdot k   }  \widehat f_0 \, dv 
    \longrightarrow \int_{\RR^N}  \widehat f_0 \, dv = \widehat \rho_0
$$
which proves (\ref{eq:lF0}).

\subsection{Passing to the limit in $a^\eps$}


The main goal of this section is the proof of the following
proposition:
\begin{lemma}\label{lem:a}
  Recall that $\gamma<2$ and $\theta(\eps)=\eps^\gamma\vphi(\eps)$,
  and let
  $$ a^\eps(p,k):= \frac{1}{\theta(\eps)} 
     \int_{\RR^N} \left[ \frac{ \nu(v)}{ \nu+\theta(\eps)  \,  p +
         \eps \, iv\cdot k   } - 1\right] \nu(v) F(v) \, dv .
  $$
  Then
  $$
  a^\eps(p,k)\longrightarrow -p- \kappa|k|^{\gamma}
  $$
  when $\eps$ goes to zero, locally uniformly with respect to $p\geq 0$
  and $k\in\RR^N$, with $\kappa$ given by (\ref{eq:D4}).  Furthermore,
  $a^\eps$ is locally bounded in $[0,\infty]\times\RR^N$ uniformly
  w.r.t. $\eps$: There exists a constant $C$ such that
  $$ |a^\eps(p,k)| \leq |p|+ C(1+ |k|^2).$$
\end{lemma}

\medskip
\noindent{\bf Proof of Lemma \ref{lem:a}.}
Observing that 
\begin{eqnarray}\label{bdd:nu-1}
\left|  \frac{ \nu}{\nu+ \theta(\eps)   p + \eps iv\cdot k  } - 1\right| \leq 2,
\end{eqnarray}
we see that $a^\eps$ is well defined for any $\eps>0$. Next, we write
\begin{eqnarray}
1 - \frac{ \nu}{\nu+ \theta(\eps)   p + \eps iv\cdot k  }  
&=&    \frac{ \nu +\theta(\eps)   p}{ (\nu+ \theta(\eps)   p)^2 
     + (\eps v\cdot k  )^2 }  \, \theta(\eps)   p\nonumber  \\
& &+ \,\,\,  \frac{ \eps^2 (v\cdot k)^2}{ (\nu+ \theta(\eps)   p)^2 
       + (\eps v\cdot k  )^2 } \nonumber \\
& &+  \,\,\,  \frac{ \eps\nu  i v\cdot k}{ (\nu+ \theta(\eps)   p)^2 
             + (\eps v\cdot k  )^2 }  . \label{eq:decomp}
\end{eqnarray}
Using the fact that $F(-v)=F(v)$ and $\nu(-v)=\nu(v)$, we deduce
\begin{eqnarray}
a^\eps(p,k ) \nonumber
&=& \label{eq:aeps}
- p \int_{\RR^N} \frac{ \nu(v) +\theta(\eps)   p}{ (\nu(v)+
  \theta(\eps)   p)^2 
     + \eps^2 (v\cdot k  )^2 }  \nu (v)F(v)  \, dv \\
&& 
- \frac{1}{\theta(\eps)}  \int_{\RR^N} \frac{(\eps v\cdot k)^2}{
  (\nu(v)+  \theta(\eps)   p)^2 
   + (\eps v\cdot k  )^2 }\nu(v)  F(v) \, dv.
\end{eqnarray}
The dominated convergence theorem immediately implies that the first
term in the right hand side converges to $-p \int_{\RR^N} F(v)\, dv =
-p$ as $\eps$ goes to zero. Furthermore, that term is clearly bounded
(in absolute value) by $|p|$.  So it only remains to show that (and
here we recall that $\theta(\eps) := \eps^{\gamma}\vphi(\eps)$)
$$
d^\eps(p,k) := \frac{1}{\eps^{\gamma}\vphi(\eps)}  
 \int_{\RR^N} \frac{(\eps\,v\cdot k)^2}{ (\nu+  \eps^{\gamma}
   \vphi(\eps)   p)^2 
 + (\eps v\cdot k  )^2 }\nu  F \, dv 
$$
converges to $ \kappa|k|^\gamma$ and is locally bounded when $\eps$
goes to zero.

For some $M>0$, we write 
$$ d^\eps(p,k)=d^\eps_1(p,k)+d^\eps_2(p,k)$$
where
\begin{eqnarray*}
d^\eps_1(p,k)& = & \frac{1}{\eps^{\gamma}\vphi(\eps)}  \int_{|v|\leq M} \frac{(\eps v\cdot k)^2}{ (\nu(v)+  \eps^{\gamma} \vphi(\eps)   p)^2 + (\eps v\cdot k  )^2 }\nu(v)  F(v) \, dv\\
& \leq & \eps^{2-\gamma}\vphi(\eps)^{-1}  \int_{|v|\leq M} \nu(v)^{-1} |v\cdot k|^2  F(v) \, dv\\
&\leq & |k|^2 \eps ^{2-\gamma}\vphi(\eps)^{-1} 
    \int_{|v|\leq M}|v|^2 \, \nu(v)^{-1} \, F(v) \, dv\\
&\leq & C \, |k|^2 \, \eps ^{2-\gamma} \, \vphi(\eps)^{-1}
\end{eqnarray*}
(using the assumption (A2) that $|v|^2 \, \nu(v)^{-1} \, F$ is
locally integrable) and
$$d_2^\eps(p,k)=\frac{1}{\eps^{\gamma} \, \vphi(\eps)}  
\int_{|v|\geq M} \frac{(\eps \, v\cdot k)^2}{ (\nu(v)+  \eps^{\gamma} \, 
  \vphi(\eps)   p)^2 
+ (\eps \, v\cdot k  )^2 }\nu(v)  F(v) \, dv.$$
\medskip

It is readily seen (using (\ref{eq:s}) and the fact that $\gamma<2$)
that
$$ d_1^\eps(p,k) \longrightarrow 0 \mbox{ as } \eps\rightarrow 0.$$
Furthermore, $d_1^\eps$ is bounded, for $\eps$ small enough by $C
|k|^2$.  So we only need to evaluate the limit of $d_2^\eps(p,k)$.
For that purpose, we first rewrite $d_2^\eps$ as follows:
$$d_2^\eps(p,k)=\frac{1}{\eps^{\gamma}\vphi(\eps)}  
\int_{|v|\geq M} \frac{ (|v|^{-\beta} \, \eps \, v\cdot
  k)^2}{ (\widetilde \nu(v)+ |v|^{-\beta} \, \eps^{\gamma} \, \vphi(\eps)
  \, p)^2 
+ (|v|^{-\beta} \, \eps \, v\cdot k  )^2 } \frac{\widetilde \nu(v)  
 \, \widetilde F_0(v) \, \ell(|v|)}{|v|^{N+\alpha -\beta}} \, dv.$$
where
$$ 
\widetilde \nu(v) = |v|^{-\beta} \, \nu(v) \quad  \mbox{ and } 
\quad \widetilde F_0 (v) = |v|^{N+\alpha} \, F_0 (v).
$$
Note that Assumptions (B1-B2-B3) (see also (\ref{eq:nub}))
imply that $\widetilde\nu$ and $\widetilde F_0$ are uniformly bounded
from above and below for $|v|\geq M$ for a suitable $M >0$, and that
\begin{equation}\label{eq:nf} 
\lim_{|v|\rightarrow \infty} \widetilde \nu (v)= \nu_0 
\quad \mbox{ and }\quad \lim_{|v|\rightarrow \infty} \widetilde F_0 (v)= \kappa_0.
\end{equation}

We now do the change of variable 
$$ w = \eps |k| |v|^{-\beta} v,$$
for which we have:
$$ |v| = \frac{|w|^{\frac{1}{1-\beta}} } { (\eps
  |k|)^{\frac{1}{1-\beta}}, }\qquad v =
\frac{w}{|w|^{\frac{-\beta}{1-\beta}}
  (\eps |k|)^{\frac{1}{1-\beta}}}, \quad  \mbox{ and
} \quad dv= \frac{(1-\beta)^{-1}}{
  (\eps|k|)^{\frac{N}{1-\beta}}
  |w|^{\frac{-\beta N}{1-\beta}}}\, dw.$$
We obtain (with $e=k/|k|$):
\begin{eqnarray*}
\!\!\!\!\!\!\!\!\!\!\!\!&& \!\!\!\!\!\!\!\!\!\!d_2^\eps(p,k) = \\
& &\!\!\!\!= \frac{(1-\beta)^{-1} }{\eps^{\gamma}\vphi(\eps)}  \int_{|w|\geq M^{1-\beta} \eps|k|} \;\; \frac{ (w\cdot e )^2}{ (\widetilde \nu^\eps (w) + |w|^{\frac{{-\beta}}{1-\beta}}(\eps |k|)^{\frac{\beta}{1-\beta}} \eps^{\gamma} \vphi(\eps)   p)^2 + (w\cdot e  )^2 } \\[4pt]
&& \qquad \times \frac{\widetilde \nu^\eps(w)  \widetilde F^\eps_0(w) } {|w|^\frac{N+\alpha-\beta}{1-\beta}} (\eps|k|)^\frac{N+\alpha-\beta}{1-\beta}   \ell\left(\frac{|w|^\frac{1}{1-\beta}} {(\eps|k|)^\frac{1}{1-\beta}} \right)\,  
(\eps|k|)^{-\frac{N}{1-\beta}} |w|^{-\frac{N-\beta}{1-\beta}}\, dw.
\end{eqnarray*}
where
$$ \widetilde \nu^\eps (w) =\widetilde \nu \left( \frac{w}{|w|^{\frac{-\beta}{1-\beta}} (\eps |k|)^{\frac{1}{1-\beta}}}\right)$$
and 
$$ \widetilde F_0^\eps (w) =\widetilde F_0 \left( \frac{w}{|w|^{\frac{-\beta}{1-\beta}} (\eps |k|)^{\frac{1}{1-\beta}}}\right)$$
(we dropped the dependence in $k$ in $  \widetilde \nu^\eps$ and $ \widetilde F_0^\eps$ to keep notation simpler).
The definition of $\gamma$ thus gives:
\begin{eqnarray}
&&\!\!\!\!\!\!\!\!\!\! d_2^\eps(p,k)=\nonumber \\
 & & \!\!\!\!= (1-\beta)^{-1} |k|^\gamma  \int_{|w|\geq M^{1-\beta}\eps|k|} \frac{ (w\cdot e )^2}{ (\widetilde \nu^\eps (w) + |w|^{\frac{-\beta}{1-\beta}}|k|^{\frac{\beta}{1-\beta}} \vphi(\eps) \eps^{\frac{\alpha}{1-\beta}} p)^2 + (w\cdot e  )^2 } \nonumber\\
&& \quad\times \frac{\widetilde \nu^\eps(w)  \widetilde F^\eps_0(w) } {|w|^{N+\gamma}}  \ell\left(\frac{|w|^\frac{1}{1-\beta}} {(\eps|k|)^\frac{1}{1-\beta}} \right)\frac{1}{\vphi(\eps)} \, dw.\label{eq:d22}
\end{eqnarray}

Finally, we note that for any $w$ and $k$, we have (note that $\beta<1$)
$$ \lim_{\eps \rightarrow 0}  \widetilde \nu^\eps (w) = \nu_0$$
$$ \lim_{\eps \rightarrow 0}  \widetilde F_0^\eps (w) = \kappa_0$$
and 
$$ \lim_{\eps \rightarrow 0} \frac{1}{\vphi(\eps)}
\ell\left(\frac{|w|^\frac{1}{1-\beta}} {(\eps|k|)^\frac{1}{1-\beta}} \right)=\lim_{\eps \rightarrow 0} \frac{1}{\ell(\eps^{-\frac{1}{1-\beta}})}
\ell\left(\frac{|w|^\frac{1}{1-\beta}} {(\eps|k|)^\frac{1}{1-\beta}} \right)=1.$$
Thus, the integrand in (\ref{eq:d22})
converges pointwise to 
$$  \frac{ (w\cdot e )^2}{ \nu_0 ^2 + (w\cdot e  )^2 }  \frac{\nu_0   \kappa_0 } {|w|^{N+\gamma}}  .
$$
Furthermore, it is bounded by
$$ C  \frac{ (w\cdot e )^2}{ \nu_1 ^2 + (w\cdot e  )^2 } \frac{1}{|w|^{N+\gamma}}  \left(1+ \left(\frac{|w|}{|k|}\right)^{\delta} \right) $$
thanks to the following lemma and an appropriate choice of $M$, see 
the following lemma (extracted from \cite[Theorem~1.5.6]{BGT}):

\begin{lemma}[Potter's Theorem]
\label{lem:potter}
  For all $\delta>0$ there exists some constants $C,M>0$ such that
$$ \left| \frac{\ell(\lambda s)}{\ell(s)} \right| \leq C (1+|\lambda
|^{\delta} ) 
\qquad \mbox{ for all  } s\geq  M \mbox{ and } \lambda \geq M/s.$$
Moreover for any $\zeta>0$:
\begin{equation}\label{eq:s}
  s^{\zeta}  \ell (s)  \longrightarrow \infty \quad \mbox{ and }
  \quad 
  s^{-\zeta}  \ell (s)  \longrightarrow 0\quad  \mbox{ as } \quad s \rightarrow \infty.
\end{equation}
\end{lemma}
Choosing $0< \delta \leq \gamma/2$ (and thus fixing a corresponding
constant $M$), we deduce that
$$| d_2^\eps(p,k)|  \leq C \, |k|^\gamma \, (1+|k|^{-\delta})$$
and by Lebesgue's dominated convergence theorem we get:
\begin{eqnarray*}
d_2^\eps(p,k)\longrightarrow (1-\beta)^{-1} \, |k|^\gamma \, 
\int_{\RR^N} \frac{ (w\cdot e )^2}{ \nu_0^2 + (w\cdot e  )^2 } 
\frac{\nu_0  \kappa_0}{ |w|^{N+\gamma} }\, dw = \kappa \, |k|^\gamma,
\end{eqnarray*}
which concludes the proof.\qed

\subsection{The remaining term}\label{sec:remaining}

We will also need the following technical lemma to pass to the limit
in the first and last term in (\ref{eq:pp}). Note that we directly
include in this lemma the limit case $\gamma =2$ since we shall need
it this same technical lemma in this limit case for the proof of
Theorem~\ref{thm:2}. 

\begin{lemma}\label{lem:aa} 
  Assume (A1-A2) and (B1-B2-B3) with $\alpha > 0$, $\beta< \alpha$ and
  $\beta \le 2-\alpha$. 
  
  Then
  \begin{eqnarray}\label{eq:lg0} 
    & &\!\!\!\!\!\!\!\! \!\!\!\!\!\!\!\!  \int_{\RR^N} 
            \left| \frac{ \nu(v)}{ \nu+  \theta(\eps)  \,  p 
                + \eps \, iv\cdot k   } - 1\right| 
    \nu(v) F(v) \, dv \le C(p,k) \, \eps^{\min\{\gamma/2 + \eta \, ;
      \, 1\}} 
  \end{eqnarray}
  for some $\eta > 0$, all $k\in\RR^N,\, p\in\RR_+$, and $C(p,k) \in
  L^\infty _{loc} (\RR_+ \times \RR^N)$.
\end{lemma}

\begin{remark}
  Note that the assumptions of Lemma \ref{lem:aa} imply that $0<
  \gamma \le 2$, and therefore $\gamma =2$ is included in these
  assumptions, except when $\alpha =1$. In other words we exclude the
  case $\alpha=\beta=1$ and hence when $\gamma=2$ we always suppose
  $\alpha >1$.
\end{remark}

\noindent{\bf Proof of Lemma~\ref{lem:aa}. }

\smallskip To prove (\ref{eq:lg0}), we split the integral into three part
parts:
 $$
\int_{\RR^N} \left| \frac{ \nu(v)}{ \nu +  \theta(\eps)  \,  p 
    + \eps \, iv\cdot k   } - 1\right| 
\nu(v) F(v) \, dv=  c_+^\eps(p,k) + c_-^\eps(p,k)$$
where
$$
c_\pm^\eps(p,k):= \int_{\left\{ |v|^{\pm 1}
  \leq \; \eps^{\frac{\mp 1}{1- \beta}}\right\}}  \left| \frac{
  \nu(v)}{ \nu+  \theta(\eps)  
\,  p + \eps \, iv\cdot k   } - 1\right| 
\nu(v) F(v) \, dv
$$
(note that this splitting is well-defined since $\beta \neq 1$ in the assumptions).

In order to get a bound for $c_+^\eps$, we note that (\ref{eq:decomp}) implies:
\begin{eqnarray*}
 \left| \frac{ \nu}{\nu+ \theta(\eps)   p + \eps \, i \, v\cdot k  } - 1\right| & \leq  & 
\nu^{-1}  \, \theta(\eps) \, |p|    + \eps^2 \nu^{-2}  (v\cdot k)^2+ \eps
\, \nu^{-1} \, | v\cdot k|
\end{eqnarray*}
and thus 
\begin{eqnarray*}
c_+^\eps(p,k) 
& \leq &  \theta(\eps) \, |p| \, \int_{\RR^N}  F(v)\, dv \\
& & 
+(\eps \, |k|)^2 \, \int_{\left\{|v|\leq \,
    \eps^{\frac{-1}{1-\beta}}\right\}}  |v|^2 \, \nu^{-1}(v) \, F(v)\, dv \\
&& + \eps \, |k| \, \int_{\left\{|v|\leq \; \eps^{\frac{-1}{1-\beta}}\right\}} |v| \,  F(v)\, dv.
\end{eqnarray*}
We write
\begin{eqnarray*}
c_+^\eps(p,k) \!&  \leq & \!  \theta (\eps) \, |p| + (\eps \, |k|)^2 \, c_{+,2}^\eps(p,k) 
+  \eps \, |k| \, c_{+,3}^\eps(p,k).
\end{eqnarray*}
The computation of $c_{+,2} ^\eps(p,k)$ requires the computation of
the following integral:
$$
\int_{\left\{|v|\leq \;\eps^{\frac{-1}{1-\beta}}\right\}}  |v|^2 \, \nu^{-1}(v)
\, F(v)\, dv
$$
which we estimate using the assumption (A2) that $|v|^2 \, \nu(v)^{-1}
\, F$ is locally integrable, and the fact that $F (v) \le C_1 \,
\ell(v) \, |v|^{-N-\alpha}$ and $\nu(v) \ge C_2 \, |v|^\beta$ for
$|v|\ge M$. When $\beta \neq 2 - \alpha$, it yields
$$
c_{+,2}^\eps \le  C \,  \left(1+\eps^{\gamma - 2 - \delta}\right)
$$
for some $\delta$ as small as we want, using the Potter's
lemma~\ref{lem:potter} to estimate $\ell$.
In the case when $\beta = 2 - \alpha$, it yields
$$
c_{+,2}^\eps \le C  \,  \left(1+ \Xi\left(\eps^{-1/(1-\beta)}\right) \right)
$$
with 
$$
\Xi(z) := \int_M ^z \frac{\ell(r)}{r} \, dr  \,\, \hbox{if} \,\, z \ge
M, \quad \Xi(z) := 0  \,\, \hbox{if} \,\, z \le M.
$$
In both cases this yields, by fixing $\delta = \gamma/2$, using the
trivial estimate 
$$
z^{-\delta} \, \Xi(z) \xrightarrow[z \to +\infty]{} 0
$$
for any $\delta >0$, and keeping only higher order terms, we obtain
the following (non optimal) bound:
$$
(\eps \, |k|)^2 \, c_{+,2}^\eps(p,k) \le C \, \eps^{\gamma /2 + \eta}
$$
for some $\eta >0$.

Similarly we found
$$
c_{+,3}^\eps  \le  C \,  \left(1+\eps^{\gamma - 1-\delta}\right) \,\, 
\hbox{if} \,\, \alpha \neq 1,
\quad c_{+,3}^\eps \le 
C  \,  \left(1+ \Xi\left(\eps^{-1/(1-\beta)}\right) \right) \,\, \hbox{if} \,\, \alpha=1
$$
for some $\delta$ as small as we want.

Using the same arguments as before (and recalling that when
$\alpha=1$, $\gamma=2$ is not allowed since $\beta < \alpha$), we get
$$
(\eps \, |k|) \, c_{+,3}^\eps(p,k) \le C \, (\eps |k|)^{\min\{ \gamma
  /2 + \eta \; ; \; 1\}}
$$
for some $\eta >0$.

\smallskip

Finally, in order to get a bound on $c_-^\eps$, we do the change of
variable $w = \eps^{\frac{1}{1-\beta}} \, v$ which leads to
\begin{eqnarray*}
&& \!\!\!\!\!\!\!\!\!c_-^\eps(p,k) \le 
 \int_{\left\{|v|\geq \; \eps^{\frac{-1}{1-\beta}}\right\}}  2 \, 
\nu(v) F(v) \, dv\\
& \le&  2 \, \eps^{\frac{\alpha-\beta}{1-\beta}}
 \int_{|w|\geq 1} 
\frac{\widetilde {\nu^\eps}(w) \widetilde F_0^\eps(w)
}{|w|^{N+\alpha-\beta}} 
    {\ell(|w|\, \eps^\frac{-1}{1-\beta})}\, dw\\
& &\!\!\!\!\!\!\leq C \, \eps^{\gamma} \, \vphi(\eps)
 \int_{|w|\geq 1} |w|^{-N-\alpha+\beta}
 \left(1+|w|^\delta\right) \, dw \\
& &\!\!\!\!\!\!\leq C' \, \eps^{\gamma} \, \vphi(\eps)
\end{eqnarray*}
where we used Potter's lemma \ref{lem:potter} for some $\delta >0$.

Collecting all the terms, and keeping only higher order terms, we find
easily (\ref{eq:lg0}).  \qed

\medskip

Finally we show that the last term in (\ref{eq:pp}) goes to zero in
$\mathcal D'([0,\infty]\times\RR^N)$: First, we note that for any
$g\in L^2(\nu F^{-1})$, we have by Cauchy-Schwarz inequality and
assumption (\ref{bdd:bnuF})
\begin{eqnarray*} 
|K(g)| & \leq &  \int_{\RR^N} \sigma(v,v') |g(v')|\, dv' \\
& \leq &  \left(\int_{\RR^N} \sigma^2 \, {F' \over \nu'}  \, dv'\right) ^{1/2}
\left(\int_{\RR^N}\frac{|g(v')|^2}{F(v')} \, \nu' \, dv'\right)^{1/2}\\
&  \leq & MÊ\,  \nu \,  F \, \|g\|_{L^2(\nu F^{-1})}.
\end{eqnarray*}

Next, using Lemma \ref{lem:aa}, we deduce that
\begin{eqnarray*}
&& \left| \frac{1}{\theta(\eps)   }\int_{\RR^N} 
\left[ \frac{\nu(v)}{ \nu(v)+  \theta(\eps)   p 
+ \eps iv\cdot k  } -1\right]K(  \widehat {g^\eps})(v)  \, dv\right| \\
&& \qquad \leq \frac{C}{\theta(\eps)} \left( \int_{\RR^N} 
\left| \frac{\nu(v)}{ \nu(v)+  \theta(\eps)   p 
+ \eps iv\cdot k  } -1\right| \nu(v) F(v) \, dv \right) \; \|\widehat
{g^\eps}\|_{L^2(\nu F^{-1})} \\
&& \qquad\leq \frac{C(p,k)}{\theta(\eps)} \, \eps^{\gamma/2 + \eta} \; \|\widehat
{g^\eps} \|_{L^2(\nu F^{-1})}.
\end{eqnarray*}
Note that here since $\gamma<2$, one can find $\eta>0$ such that
$\gamma/2+\eta <1$, and therefore the estimate (\ref{eq:lg0}) in Lemma
\ref{lem:aa} yields a bound from the below of the form $C \, \eps^{\gamma/2+\eta}$.

Last, recalling inequality (\ref{eq:gg}), we have for any test function 
$\phi\in \mathcal D((0,\infty)\times\RR^N)$:
\begin{eqnarray*}
&& \left| \int_0^\infty\int_{\RR^N }\phi(p,k) \frac{1}{\theta(\eps)} 
\int_{\RR^N} \left[ \frac{\nu(v)}{ \nu(v)+  \theta(\eps)   p + \eps
    iv\cdot k  } 
    -1\right] \, K(  \widehat {g^\eps})(v)  \, dv \, dk\, dp \right|\\
&& \qquad\leq \frac{C \, \eps^{\gamma/2 + \eta}}{\theta(\eps)} \,
\left(  \int_0^\infty  \left( \int_{\RR^N } C(p,k)^2 \, 
            \phi(p,k)^2 \, dk\,\right)^{1/2}  dp \right) \\[5pt]
&&\qquad \quad \quad \quad \times \|\widehat {g^\eps} \|_{L^\infty((a,\infty);L^2(\nu F^{-1}))} 
 \\
&& \qquad\leq \frac{C(\phi) \, \eps^{\gamma/2 + \eta}}{\theta(\eps)}
\,  \|{g^\eps} \|_{L^2((0,\infty); L^2(\nu F^{-1}))} 
 \\
&& \qquad\leq C(\phi) \, \eps^{\eta} \, \varphi(\eps)^{-1} 
\end{eqnarray*}
which goes to zero thanks to (\ref{eq:s}) and the fact that  $\eta>0$.
\qed

\subsection{Conclusion}
We now have all the tools to pass to the limit $\eps\rightarrow 0$ in
(\ref{eq:pp}). The first part of Lemma \ref{lem:aa} implies that the
first term converges to $\widehat{\rho_0}(k)$ locally uniformly with
respect to $k$ and $p$, while the second term can be rewritten as
$a^\eps(p,k) \widehat {\rho^\eps}$ which converges, using Lemma
\ref{lem:a}, in $\mathcal D'((0,\infty)\times\RR^N)$ to $(-p-
\kappa|k|^\gamma) \widehat \rho$.  

The remaining term goes to $0$ thanks the the second part of Lemma
\ref{lem:aa} and the preceeding subsection. 

Hence we can pass to the limit in (\ref{eq:pp}) and recover
(\ref{eq:laplace-fourier}). The proof of Theorem \ref{thm:1} can then
be completed as the proof of Theorem~\ref{thm:0}.  


\section{Proof of Theorem~\ref{thm:2}}\label{sec:5}

The proof of Theorem \ref{thm:3} is very similar to the proof of
Theorem~\ref{thm:1}. We recall the main steps for the reader's sake.

We recall that we define the following time scale in the critical case:
$$
\theta(\eps) = \eps^2 \, \ell(\eps^{-\frac{1}{1-\beta}}) \, \ln(\eps^{-1}).
$$
and we shall use the same notation as in the previous section:
$$
\varphi(\eps) =  \ell(\eps^{-\frac{1}{1-\beta}}).
$$
\medskip

The starting point is again the equality (\ref{eq:pp}). As in the
previous section, we see that if we take that the first term converges
to $\rho_0$ when $\eps$ goes to zero.

Proceeding as in Section \ref{sec:remaining} we can also
show (using Lemma \ref{lem:aa} with $\gamma=2$) that the last term
goes to zero:   
\begin{eqnarray*}
&& \left| \frac{1}{\theta(\eps)   }\int_{\RR^N} 
\left[ \frac{\nu(v)}{ \nu(v)+  \theta(\eps)   p 
+ \eps iv\cdot k  } -1\right]K(  \widehat {g^\eps})(v)  \, dv\right| \\
&& \qquad \leq \frac{C}{\theta(\eps)} \left( \int_{\RR^N} 
\left| \frac{\nu(v)}{ \nu(v)+  \theta(\eps)   p 
+ \eps iv\cdot k  } -1\right| \nu(v) F(v) \, dv \right) \; \|\widehat
{g^\eps}\|_{L^2(\nu F^{-1})} \\
&& \qquad\leq \frac{C(p,k)}{\theta(\eps)} \, \eps \; \|\widehat
{g^\eps} \|_{L^2(\nu F^{-1})}.
\end{eqnarray*}
Note that here since $\gamma=2$, the estimate (\ref{eq:lg0}) in Lemma
\ref{lem:aa} yields a bound from the below of the form $C \,\eps$.

Last, recalling inequality (\ref{eq:gg}), we have for any test function 
$\phi\in \mathcal D((0,\infty)\times\RR^N)$:
\begin{eqnarray*}
&& \left| \int_0^\infty\int_{\RR^N }\phi(p,k) \frac{1}{\theta(\eps)} 
\int_{\RR^N} \left[ \frac{\nu(v)}{ \nu(v)+  \theta(\eps)   p + \eps
    iv\cdot k  } 
    -1\right] \, K(  \widehat {g^\eps})(v)  \, dv \, dk\, dp \right|\\
&& \qquad\leq \frac{C \, \eps}{\theta(\eps)} \,
 \left(  \int_0^\infty  \left( \int_{\RR^N } C(p,k)^2 \, 
            \phi(p,k)^2 \, dk\,\right)^{1/2}  dp \right) 
 \, \|\widehat {g^\eps} \|_{L^\infty((a,\infty);L^2(\nu F^{-1}))} 
 \\
&& \qquad\leq \frac{C(\phi) \, \eps}{\theta(\eps)} \,  \|{g^\eps} \|_{L^2((0,\infty); L^2(\nu F^{-1}))} 
 \\
&& \qquad\leq \frac{C(\phi)}{\varphi(\eps) \, \ln(\eps^{-1})} 
\end{eqnarray*}
which goes to zero thanks to assumption (\ref{ellcritic}) in Theorem
\ref{thm:2}: $\ell(r) \, \ln r \to +\infty$ as $r \to +\infty$ indeed
immediately implies that $\ell(\eps^{-1/(1-\beta)}) \, \ln(\eps^{-1})
\to + \infty$ as $\eps \to 0$ (recall that here $1 - \beta = \alpha -
1 >0$).

Finally, we are left with the task of determining the
limit of the symbol $a^\eps$ in (\ref{eq:aeps}).

It has already been proved that the first term in the right hand side
of (\ref{eq:aeps}) is bounded (uniformly in $\eps$) by $|p|$ and
converges to $-p \int_{\RR^N} F(v)\, dv = -p$ as $\eps$ goes to zero.

So it only remains to show that 
$$
d^\eps(p,k) := \frac{1}{\theta(\eps)} \,
\int_{\RR^N} \frac{(\eps\,v\cdot k)^2}{ (\nu+  \theta(\eps) p)^2 
+ (\eps v\cdot k  )^2 } \, \nu(v) \, F(v) \, dv 
$$
converges to $\kappa \, |k|^2$ and is locally bounded when $\eps$ goes to
zero.

For some $M>0$, we write 
$$ d^\eps(p,k)=d^\eps_1(p,k)+d^\eps_2(p,k),$$
where
\begin{eqnarray*}
d^\eps_1(p,k)& = & \frac{1}{\theta(\eps)} \, 
 \int_{|v|\leq M} \frac{(\eps v\cdot k)^2}{ (\nu(v)+  \theta(\eps) \, p)^2 
 + (\eps v\cdot k  )^2 } \, \nu(v) \,  F(v) \, dv\\
& \leq & \frac{1}{\varphi(\eps) \, \ln(\eps^{-1})}  \, \int_{|v|\leq M}\nu(v)^{-1} \, |v\cdot k|^2  F(v) \, dv\\
&\leq & \frac{|k|^2}{\varphi(\eps) \, \ln(\eps^{-1})} \,
 \int_{|v|\leq M}|v|^2 \, \nu(v)^{-1}  F(v) \, dv\\
&\leq & \frac{C \, |k|^2}{\varphi(\eps) \, \ln(\eps^{-1})} 
\end{eqnarray*}
(using the assumption (A2) that $|v|^2 \, \nu(v)^{-1} \, F$ is locally
integrable) and 
$$
d_2^\eps(p,k)=\frac{1}{\theta(\eps)} \,  
\int_{|v|\geq M} \frac{(\eps \, v\cdot k)^2}{ (\nu(v)+  \theta(\eps)
  \, p)^2  + (\eps \, v\cdot k  )^2 } \, \nu(v) \, F(v) \, dv.
$$
\medskip

Concerning $d_1^\eps(p,k)$, as above we deduce from the assumption
(\ref{ellcritic}) on $\ell$ that 
$$
\varphi(\eps) \, \ln(\eps^{-1}) \to +\infty \quad \mbox{ as
} \eps \to 0,
$$
and therefore $d_1^\eps(p,k)$ goes to $0$ as $\eps$ goes to $0$.
Moreover it is clearly bounded, for $\eps$ small enough, by some $C \,
|k|^2$.

Now it remains to evaluate the limit of $d_2^\eps(p,k)$.  For that
purpose, we first rewrite $d_2^\eps$ as follows:
$$
d_2^\eps(p,k)=\frac{1}{\theta(\eps)}  
\int_{|v|\geq M} \frac{ (|v|^{-\beta} \eps \, v\cdot k)^2}
{ (\widetilde \nu(v)+ |v|^{-\beta} \, \theta(\eps) \, p)^2 
+ (|v|^{-\beta} \, \eps \, v \cdot k  )^2 } 
  \frac{\widetilde \nu(v)  \, \widetilde
    F_0(v)}{|v|^{N+\alpha -\beta}} 
   \, \ell(|v|) \, dv.$$
where
$$ 
\widetilde \nu(v) = |v|^{-\beta} \, \nu(v) \quad  
\mbox{ and } \quad \widetilde F_0 (v) = |v|^{N+\alpha} \, F_0 (v).
$$
Note that Assumptions (B1-B2-B3) (see also (\ref{eq:nub})) imply
that $\widetilde\nu$ and $\widetilde F_0$ are uniformly bounded from
above and below for $|v|\geq M$ and that
$$
\lim_{|v|\rightarrow \infty} \widetilde \nu (v) 
= \nu_0 \quad \mbox{ and }\quad \lim_{|v|\rightarrow \infty} \widetilde F_0 (v)= \kappa_0.
$$
We do again the change of variable $ w = \eps |k|
|v|^{-\beta} v$ (recall again that $\beta \neq 1$
from the assumptions), and we obtain (with $e=k/|k|$):
\begin{eqnarray*}
\!\!\!\!\!\!\!\!\!\!\!\!&& \!\!\!\!\!\!\!\!\!\!d_2^\eps(p,k) = \\
& &\!\!\!\!= \frac{(1-\beta)^{-1} }{\theta(\eps)} 
\int_{|w|\geq M^{1-\beta} \eps|k|} \;\; 
\frac{ (w\cdot e )^2}{ \left(\widetilde \nu^\eps (w) 
+ |w|^{\frac{-\beta}{1-\beta}}(\eps |k|)^{\frac{\beta}{1-\beta}}
\eps^{2} \varphi(\eps) p\right)^2 + (w\cdot e  )^2 } \\[4pt]
&& \qquad \times \frac{\widetilde \nu^\eps(w)  \widetilde F^\eps_0(w)}
{|w|^\frac{N+\alpha-\beta}{1-\beta}} (\eps|k|)^\frac{N+\alpha-\beta}{1-\beta}  \,  
 \ell\left(\frac{|w|^\frac{1}{1-\beta}} {(\eps|k|)^\frac{1}{1-\beta}}
 \right) \, (\eps|k|)^{-\frac{N}{1-\beta}} \, |w|^{-\frac{N-\beta}{1-\beta}}\, dw,
\end{eqnarray*}
where
$$ \widetilde \nu^\eps (w) =\widetilde \nu \left(
  \frac{w}{|w|^{\frac{-\beta}{1-\beta}} (\eps |k|)^{\frac{1}{1-\beta}}}\right)$$
and 
$$ \widetilde F_0^\eps (w) =\widetilde F_0 \left(
  \frac{w}{|w|^{\frac{-\beta}{1-\beta}} (\eps
    |k|)^{\frac{1}{1-\beta}}}\right)
$$ 
(we dropped the dependence in $k$ in $\widetilde \nu^\eps$ and $
\widetilde F_0^\eps$ to keep notation simpler). 

The fact that $\gamma=2$ now yields
\begin{eqnarray}
  &&\!\!\!\!\!\!\!\!\!\! d_2^\eps(p,k)=\nonumber \\
  & & \!\!\!\!= \frac{(1-\beta)^{-1} }{\varphi(\eps) \, \ln(\eps^{-1})}|k|^2
  \int_{|w|\geq M^{1-\beta}\eps|k|} \frac{ (w\cdot e )^2}{ \left(\widetilde
    \nu^\eps (w) +
    |w|^{\frac{-\beta}{1-\beta}}|k|^{\frac{\beta}{1-\beta}}
    \varphi(\eps) \eps^{\frac{\alpha}{1-\beta}} p\right)^2 + (w\cdot e  )^2 } \nonumber \\
  && \quad\times \frac{\widetilde \nu^\eps(w)  \widetilde F^\eps_0(w)}
  {|w|^{N+2}} \,  \ell\left(\frac{|w|^\frac{1}{1-\beta}}
    {(\eps|k|)^\frac{1}{1-\beta}} \right) \, dw.\label{eq:d222}
\end{eqnarray}

Lemma \ref{lem:fv} below implies
\begin{equation*}
d^\eps_2(p,k) \xrightarrow[\eps \to 0]{}  \kappa \, |k|^2
\end{equation*}
which concludes the proof of Theorem~\ref{thm:2}.

\begin{lemma}\label{lem:fv}
We have: 
\begin{eqnarray}
  & &  \frac{(1-\beta)^{-1} }{\varphi(\eps) \, \ln(\eps^{-1})} \, 
  \int_{|w|\geq M^{1-\beta}\eps|k|} \frac{ (w\cdot e )^2}{ \left(\widetilde
    \nu^\eps (w) +
    |w|^{\frac{-\beta}{1-\beta}}|k|^{\frac{\beta}{1-\beta}}
    \varphi(\eps) \eps^{\frac{\alpha}{1-\beta}} p\right)^2 + (w\cdot e  )^2 } \nonumber \\
  && \quad\times \frac{\widetilde \nu^\eps(w)  \widetilde F^\eps_0(w)}
  {|w|^{N+2}} \,  \ell\left(\frac{|w|^\frac{1}{1-\beta}}
    {(\eps|k|)^\frac{1}{1-\beta}} \right) \, dw \to \kappa \label{eq:lemchiant}
\end{eqnarray}
as $\eps \rightarrow 0$.
\end{lemma}
{\bf Proof of Lemma \ref{lem:fv}:}
First, we show that 
\begin{equation}\label{eq:hop}
\kappa := \frac{1}{(1-\beta)} \, \lim_{\lambda\to 0} \frac{1}{\ln(\lambda^{-1})} 
\int_{|w| \ge \lambda }  \ \frac{ w_1^2}{ \nu_0 ^2 + w_1^2 }  \,
{\kappa_0 \, \nu _0 \over  |w|^{N+2}} \,dw
\end{equation}
is well defined:
We denote
$$\psi(\lambda)= \int_{|w| \ge \lambda }  \ \frac{w_1^2}{
  \nu_0 ^2 + w_1^2 }  \, {\kappa_0 \, \nu_0 \over  |w|^{N+2}} \,dw $$
and then we have
\begin{eqnarray*}
\psi'(\lambda) & =&  - \int_{|w| = \lambda }  \ \frac{ w_1^2}{ \nu_0
  ^2 + w_1^2 }  \, {\kappa_0 \, \nu_0 \over  |w|^{N+2}} \,d\sigma(w)\\
& =&  - \frac{{\kappa_0 \, \nu_0} }{\lambda^{N+2}} \int_{|w| =
  \lambda}  \ \frac{ w_1^2}{ \nu_0 ^2 + w_1^2 } \,d\sigma(w)\\
& =&  - \frac{{\kappa_0 \, \nu_0} }{\lambda} \int_{|v| = 1 }  
  \ \frac{ v_1^2}{ \nu_0 ^2 + \lambda^2 v_1^2 }  \, d\sigma(v).\\
\end{eqnarray*}
We deduce from L'Hopital's rule that
$$\kappa=\frac{1}{(1-\beta)} \,  
\lim_{\lambda\to 0} \frac{\psi(\lambda)}{\ln(\lambda^{-1})} = 
\frac{1}{(1-\beta)} \,  \lim_{\lambda\rightarrow 0}
\frac{\psi'(\lambda)}{\ln'(\lambda^{-1})} 
= {\kappa_0 \over (1- \beta) \, \nu_0} \, \int_{|v| = 1 }   v_1^2 \, d\sigma(v).$$

\smallskip

Then, we note that from the assumptions (B1-B2) and the definition of
a slowly varying function, for any $w$ and $k$, we have
$$ \widetilde \nu^\eps (w) =  \nu_0 + o(1) \quad \mbox{ as } \eps \to 0,$$
$$ \widetilde F_0^\eps (w) = \kappa_0 + o(1) \quad \mbox{ as } \eps \to 0,$$
and 
$$ \frac{1}{\vphi(\eps)}
\ell\left(\frac{|w|^\frac{1}{1-\beta}} {(\eps|k|)^\frac{1}{1-\beta}}
\right)
= \frac{1}{\ell(\eps^{-\frac{1}{1-\beta}})}
\ell\left(\frac{|w|^\frac{1}{1-\beta}} {(\eps|k|)^\frac{1}{1-\beta}}
\right)=1 + o(1) \quad \mbox{ as } \eps \to 0.$$

Thus, the integrand in (\ref{eq:d222}) divided by $\varphi(\eps)$
converges pointwise to
$$  
\frac{ (w\cdot e )^2}{ \nu_0 ^2 + (w\cdot e  )^2 }  \frac{\nu_0
  \kappa_0 } {|w|^{N+2}}  
 \quad \mbox{ as } \eps \to 0,
$$
and it is bounded by
$$ C  \frac{ (w\cdot e )^2}{ \nu_1 ^2 + (w\cdot e  )^2 } \frac{1}{|w|^{N+2}}.  $$

Then, easy but tedious computations yield, for bounded $p$ and $k$ and
for $\eps$ going to $0$:
\begin{eqnarray*}
  && \Bigg| \int_{|w|\geq M^{1-\beta}\eps|k|} \frac{ (w\cdot e )^2}{ \left(\widetilde
    \nu^\eps (w) +
    |w|^{\frac{-\beta}{1-\beta}}|k|^{\frac{\beta}{1-\beta}}
    \varphi(\eps) \eps^{\frac{\alpha}{1-\beta}} p\right)^2 + (w\cdot e  )^2 } \\
  && \qquad\times \frac{\widetilde \nu^\eps(w)  \widetilde F^\eps_0(w)}
  {|w|^{N+2}} \,  \ell\left(\frac{|w|^\frac{1}{1-\beta}}
    {(\eps|k|)^\frac{1}{1-\beta}} \right) \, dw 
  - \int_{|w| \ge  M^{1-\beta} \eps \, |k| }  \ \frac{ w_1^2}{\nu_0 ^2
    + w_1^2 }  \, {\kappa_0 \, \nu_0 \over  |w|^{N+2}} \,dw\Bigg|\\
  && \quad \leq o(1) \,  \int_{|w| \ge  M^{1-\beta} \eps \, |k| }  \ \frac{ w_1^2}{\nu_0 ^2
    + w_1^2 }  \, {\kappa_0 \, \nu_0 \over  |w|^{N+2}} \,dw \\
  && \qquad C \, \varphi(\eps) \, \eps^{\frac{\alpha}{1-\beta}} \, 
    \int_{|w| \ge  M^{1-\beta} \eps \, |k| }  \ \frac{ w_1^2}{\nu_0 ^2
    + w_1^2 }  \, {\kappa_0 \, \nu_0 \over  |w|^{N+2+\beta/(1-\beta)}}
  \,dw \\
  && \qquad \leq o(1) \, \ln(\eps^{-1}) + C \, \varphi(\eps) \,
  \eps^2 \\
  &&  \qquad \leq o(1) \, \ln(\eps^{-1}).
\end{eqnarray*}

Thus we have
\begin{eqnarray*}
&& \lim_{\eps\to 0}  \frac{1}{\ln(\eps^{-1})} \int_{|w|\geq M^{1-\beta}\eps|k|} \frac{ (w\cdot e )^2}{ \left(\widetilde
    \nu^\eps (w) +
    |w|^{\frac{-\beta}{1-\beta}}|k|^{\frac{\beta}{1-\beta}}
    \varphi(\eps) \eps^{\frac{\alpha}{1-\beta}} p\right)^2 + (w\cdot e  )^2 } \nonumber \\
  && \quad\times \frac{\widetilde \nu^\eps(w)  \widetilde F^\eps_0(w)}
  {|w|^{N+2}} \,  \ell\left(\frac{|w|^\frac{1}{1-\beta}}
    {(\eps|k|)^\frac{1}{1-\beta}} \right) \, dw  \\
&& \quad = \lim_{\eps\to 0}  \frac{1}{\ln(\eps^{-1})} \int_{|w| \ge
  M^{1-\beta} \eps |k| }  
   \ \frac{ w_1^2}{ \nu_0 ^2 + w_1^2 }  \, {\kappa_0 \, \nu_0 \over  |w|^{N+2}} \,dw .
\end{eqnarray*}

Finally, using the fact that for all $k\neq 0$,
$$\lim_{\eps\to 0} \frac{\ln((M^{1-\beta} \eps |k|)^{-1})}{\ln(\eps^{-1})} = 1,$$
we deduce (\ref{eq:lemchiant}), which completes the proof of Lemma~\ref{lem:fv}. \qed

\section{Appendix}\label{sec:A}

\begin{lemma}
\label{lem:LA} Let us assume that $b = b(v,v') = |v-v'|^\beta$ with $\beta \in \RR$.
\begin{itemize}
\item[(i)] For $\beta \in (-N,\alpha)$ there exists $C$ such
that
$$
C^{-1} \, \langle v \rangle^\beta \le \nu(v) := K(F)(v) \le C \,
\langle v \rangle^\beta, \quad \forall \, v \in \RR^N.
$$
\item[(ii)] For $\beta \in (-\min\{\alpha;N\},\min\{\alpha;N\})$ there exists $M$ such
that
$$
\int_{\RR^N} F' \, {\nu \over b} \, dv' \le M, \quad \forall \, v \in \RR^N.
$$
\item[(iii)] For $\beta \in  (-\min\{\alpha;N/2\},\alpha)$ there exists $M$ such that
$$
\int_{\RR^N} {F' \over \nu'} \, {b^2 \over \nu^2} \, dv' \le M \qquad \forall \, v \in \RR^N.
$$
\end{itemize}
\end{lemma}

\smallskip\noindent 
{\bf Proof of Lemma \ref{lem:LA}.}  {\sl Point (i). } The case $\beta
\ge 0$ being very simple, we only deal with the case $\beta < 0$ and
large values of $|v|$. We split $\nu = \nu_1 + \nu_2$ with
$$
\nu_1(v) = \int_{|v'-v| \le |v|/2} F(v') \, |v-v'|^\beta \, dv', \quad
\nu_2(v) = \int_{|v'-v| \ge |v|/2} F(v') \, |v-v'|^\beta \, dv'.
$$
On the one hand
$$
\nu_1(v) \le F(|v|/2) \int_{|v'-v| \le |v|/2} |v-v'|^\beta \, dv' \le
F_0 \, |v|^{-\alpha-N} \, |v|^{N+\beta} \le C \, |v|^{\beta}.
$$
On the other hand
$$
\nu_2(v) \le (|v|/2)^\beta \int_{|v'-v| \ge |v|/2} F' \, dv' \le C \,
|v|^{\beta}.
$$
Finally,
$$
\nu_2(v) \ge (3 \, |v|/2)^\beta \int_{|v|/2 \le |v'-v| \le 3 |v|/2} F'
\, dv' \ge (3 \, |v|/2)^\beta \int_{|v| \le 1} F' \, dv' \ge C^{-1} \,
|v|^{\beta}
$$
for $|v| \ge 2$. 

\smallskip\noindent {\sl Points (ii) and (iii). } From
Lemma~\ref{lem:LA} we have on the one hand
$$
\int_{\RR^N} {F' \over b} \, dv'= \int_{\RR^N} F' \, |v-v'|^{-\beta}
\le C \, \langle v \rangle^{-\beta} \le M \, \nu^{-1} \qquad \forall
\, v \in \RR^N,
$$
and on the other hand 
$$
\int_{\RR^N} {F' \over \nu'} \, b^2\, dv' \le \int_{\RR^N} \langle v'
\rangle^{-N-\alpha-\beta}\, |v-v'|^{2\beta} \le C \, \langle v
\rangle^{2\beta} \le M \, \nu^{2}, \quad \forall \, v \in \RR^N.
$$
\qed
\bigskip

\noindent {\bf Acknowledgements:} We thank Jean Dolbeault and Stefano
Olla for fruitful discussions during the preparation of this work.


\bigskip

\footnotesize

\begin{thebibliography}{AA}


\bibitem{BSS} (MR0743736) \newblock C. Bardos, R. Santos and R.
  Sentis, \newblock{\em Diffusion approximation and computation of the
    critical size}, \newblock Trans. A. M. S., {\bf 284} (1984),
  617--649.

\bibitem{BLP} (MR0533346) \newblock A. Bensoussan, J. L. Lions and G.
  Papanicolaou, \newblock {\em Boundary layers and homogenization of transport
    processes}, \newblock Publ. Res. Inst. Math. Sci. {\bf 15} (1979), 53--157.

\bibitem{BGT} (MR1015093) \newblock N.H. Bingham, C.M. Goldie and J.L.
  Teugels, \newblock {\em Regular variation}, \newblock Encyclopedia
  of Mathematics and its Applications {\bf 27}. Cambridge University Press,
  Cambridge, 1989.

\bibitem{BCG} (MR1749231) \newblock A.V. Bobylev, J.A. Carrillo and
  I.M. Gamba, \newblock {\em On some properties of kinetic and hydrodynamic equations
  for inelastic interactions}, \newblock J. Statist. Phys. {\bf 98}
  (2000), 743--773.

\bibitem{BG} (MR2264617) \newblock A.V. Bobylev and I.M. Gamba,
  \newblock {\em Boltzmann equations for mixtures of Maxwell gases: exact
  solutions and power like tails}, \newblock J. Stat. Phys. {\bf 124} (2006), 497--516.

\bibitem{BGTh} (MR1174046) \newblock C. B\"orgers, C. Greengard and E. Thomann,
  \newblock {\em The Diffusion Limit of Free Molecular Flow in Thin Plane
    Channels}, \newblock SIAM J. Appl. Math. {\bf 52} (1992), 1057--1075.

\bibitem{DGP} (MR1803225) \newblock P. Degond, T. Goudon and F. Poupaud,
  \newblock {\em Diffusion limit for non homogeneous and non-micro-reversibles
    processes}, \newblock Indiana Univ. Math. J. {\bf 49} (2000), 1175--1198.

\bibitem{D1} (MR1649918) \newblock C. Dogbe, \newblock {\em Diffusion
    Anormale pour le Gaz de Knudsen}. \newblock C. R. Acad. Sci. Paris
  S\'er. I Math. {\bf 326} (1998), 1025--1030.

\bibitem{D2} (MR1788479) \newblock C. Dogbe, \newblock {\em Anomalous
    diffusion limit induced on a kinetic equation}, J. Statist. Phys.
  {\bf 100} (2000), 603--632.

\bibitem{DT} \newblock D. Duering and G. Toscani, \newblock {\em
    Anomalous diffusion limit indeuced on a kinetic equation},
    Physica A {\bf 384} (2007), 493--506.

\bibitem{EB} (MR1942001) \newblock M.H. Ernst, and R. Brito, \newblock
  {\em Scaling solutions of inelastic Boltzmann equations with
    over-populated high energy tails}, \newblock J. Statist. Phys.
  {\bf 109} (2002), 407--432.

\bibitem{G} (MR1632712) \newblock F. Golse, \newblock {\em Anomalous
    Diffusion Limit for the Knudsen Gas}, \newblock Asymptotic Anal. {\bf
    17}, (1998), 1--12.

\bibitem{JKO} \newblock H. Jara, T. Komorowski and S. Olla, \newblock
  {\em A limit theorem for additive functionals of a Markov chain},
  \newblock Preprint arXiv:0809.0177.

\bibitem{LK} (MR0339741) \newblock E.W. Larsen and J.B. Keller,
  \newblock {\em Asymptotic solution of neutron transport problems for
    small mean free paths}, \newblock J. Math.  Phys. {\bf 15} (1974),
  75--81.

\bibitem{MT} \newblock A. Mellet and B. Texier, \newblock work in preparation.

\bibitem{MR} \newblock D.A. Mendis and M. Rosenberg, \newblock {\em
    Cosmic dusty plasma}, \newblock Annu. Rev. Astron. Astrophys. {\bf 32}
  (1994), 419--63.

\bibitem{N} \newblock M.E.J. Newman, \newblock {\em Power laws, Pareto
    distributions and Zipf's law}, \newblock Contemp. Phys. {\bf 46}
  (2005), 323--351.

\bibitem{ST} \newblock D. Summers and R.M. Thorne, \newblock {\em The modified plasma
    dispersion function}, \newblock Phys. Fluids {\bf 83} (1991),
  1835--1847.

\bibitem{Vill-gran} (MR2264625) \newblock C. Villani, \newblock 
  {\em Mathematics of granular materials},  \newblock
  J. Stat. Phys. {\bf 124} (2006), 781--822.

\bibitem{W} \newblock E. Wigner, \newblock {\em Nuclear reactor
    theory}, \newblock AMS (1961).

\bibitem{Wright} \newblock I. Wright, \newblock {\em The social
    architecture of capitalism}, \newblock Physica A {\bf 346} (2005),
  589--620.

\end{thebibliography}
\end{document}